\RequirePackage{etoolbox}
\csdef{input@path}{{style/}{graphics/}}
\documentclass[aap,MSNbibl,nameyear,dvips]{arximspdf}
\usepackage{graphicx}
\usepackage{url,breakurl}
\doi{10.1214/14-AAP1016} 
\volume{25}
\issue{2}
\pubyear{2015}
\firstpage{986}
\lastpage{1029}

\makeatletter
\newcommand{\eqref}[1]{(\ref{#1})}
\newcommand{\DecrAvg}[1]{{\tilde{#1}}}
\newcommand{\ApproxVar}[1]{{\tilde{#1}}}

\newcommand{\DecrAvgExtrap}[1]{{\hat{#1}}}
\newcommand{\ApproxVarExtrap}[1]{{\hat{#1}}}

\newcommand{\prob}[1][P]{{\mathbb{#1}}}
\newcommand{\tendsas}[0]{\mathop{\stackrel{a.s.}{\longrightarrow}}}
\def\mathds{\mathbb}
\newcommand{\F}{\mathcal{F}}
\newcommand{\RR}{\mathds{R}}
\newcommand{\NN}{\mathds{N}}
\newcommand{\QQ}{\mathds{Q}}

\newcommand{\tbar}[1][t]{{\underline{#1}}}

\newtheorem{MyTheorem}{Theorem}[section]
\newtheorem{MyLemma}[MyTheorem]{Lemma}
\newtheorem{MyProposition}[MyTheorem]{Proposition}
\newtheorem{MyCorollary}[MyTheorem]{Corollary}
\newproclaim{MyAssumption}[MyTheorem]{Assumption}

\newproclaim{MyDefinition}[MyTheorem]{Definition}

\newproclaim{MyExample}[MyTheorem]{Example}

\newproclaim{MyRemark}[MyTheorem]{Remark}

\newtheorem{HypSteps}{Hypothesis}
\newtheorem{HypFastScale}{Hypothesis}
\newtheorem{HypSlowScale}{Hypothesis}
\makeatother

\begin{document}
\begin{frontmatter}

\title{A decreasing step method for strongly oscillating stochastic models}
\runtitle{A decreasing step method for strongly oscillating SDEs}

\begin{aug}
\author[A]{\fnms{Camilo~Andr\'{e}s} \snm{Garc\'{i}a Trillos}\corref{}\ead[label=e1]{camilo@unice.fr}}
\runauthor{C.~A. Garc\'{i}a Trillos}
\affiliation{Universit\'{e} Nice Sophia Antipolis}

\address[A]{Universit\'{e} Nice Sophia Antipolis\\
CNRS, LJAD, UMR 7351\\
Parc Valrose\\
06108 Nice Cedex 02\\
France\\
\printead{e1}}
\end{aug}

\received{\smonth{10} \syear{2012}}
\revised{\smonth{1} \syear{2014}}

%
\begin{abstract}
We propose an algorithm for approximating the solution of a strongly
oscillating SDE, that is, a system in which some ergodic state
variables evolve quickly with respect to the other variables. The
algorithm profits from homogenization results and consists of an Euler
scheme for the slow scale variables coupled with a decreasing step
estimator for the ergodic averages of the quick variables. We prove the
strong convergence of the algorithm as well as a C.L.T. like limit
result for the normalized error distribution. In addition, we propose
an extrapolated version that has an asymptotically lower complexity and
satisfies the same properties as the original version.
\end{abstract}

%
\begin{keyword}[class=AMS]
\kwd[Primary ]{60H35}
\kwd[; secondary ]{65C30}
\end{keyword}

\begin{keyword}
\kwd{Stochastic approximation}
\kwd{strongly oscillating}
\kwd{multi-scale system}
\kwd{ergodicity}
\kwd{limit distribution}
\end{keyword}
\end{frontmatter}

\section{Introduction}\label{sec1}

Consider a system of stochastic equations of the form
%
\begin{equation}
\label{EqTheSystem} \cases{\displaystyle X_t^\varepsilon= x_0 + \int
_0^t f \bigl(X_s^\varepsilon,Y_s^\varepsilon
\bigr) \,ds + \int_0^t g \bigl(X_s^\varepsilon,Y_s^\varepsilon
\bigr) \,dW_s,\vspace*{2pt}
\cr
\displaystyle Y_t^\varepsilon=
y_0 + \varepsilon^{-1} \int_0^t
b \bigl(X_s^\varepsilon,Y_s^\varepsilon
\bigr) \,ds + \varepsilon^{-1/2} \int_0^t
\sigma \bigl(X_s^\varepsilon, Y_s^\varepsilon
\bigr) \,d\tilde{W}_s,}
\end{equation}
where $X_t^\varepsilon$ is a $d_x$-dimensional process,
$Y_t^\varepsilon$ a
$d_y$-dimensional process, $W$ and $\tilde{W}$ are two independent
Brownian motions of dimensions $d_{x}$ and $d_{y}$, and the functions
$b,\sigma, f \mbox{ and } g$ have the right dimensions.

This type of system models the dynamics of two sets of interacting
variables evolving in different time scales. The difference between
time scales is controlled by the parameter $\varepsilon$. In many domains
the most interesting case of study is that of the regime when
$\varepsilon
\ll1$, that is, the situation in which $X^\varepsilon$ is evolving very
slowly compared to $Y^\varepsilon$ (for this reason we will frequently
denominate them as \emph{slow scale} and \emph{fast scale} variables,
resp.). This regime may be studied by singular perturbation
techniques as the ones presented in \citet{bensoussanasymptotic1978}
for deterministic models: instead of looking at the system with a small
$\varepsilon$, we study the limit of \eqref{EqTheSystem} as
$\varepsilon
\rightarrow0$ (when it exists) and estimate the error induced by this
approximation.

There exist several analytical works with applications in different
domains on the described type of approximation for stochastic models.
For example in Majda, Timofeyev and
Vanden-Eijnden (\citeyear{majdamathematical2001}) a climate model is
considered and studied on the advection scale (i.e., in the time scale
of the slow variable). In \citet{fouquederivatives2000} and \citet
{fouquesingular2003} a system similar to \eqref{EqTheSystem} is
presented and studied for pricing derivatives in the context of
stochastic volatility models. A complete study with rather general
hypothesis on the coefficients of the system is found in \citet
{pardouxpoisson2001} and \citet{pardouxpoisson2003}. In these papers
a system with a fast scale ergodic diffusion is considered. More
precisely, if
%
\begin{equation}
\label{EqErgodicDiffusion} Y_t^x = y_0 + \int
_0^t b \bigl(x,Y_s^x
\bigr) \,ds + \int_0^t \sigma \bigl(x,
Y_s^x \bigr) \,d\tilde{W}_s,
\end{equation}
is ergodic with unique invariant measure $\mu^x$, we might define the
\emph{effective equation}
%
\begin{equation}
\label{EqEffectiveEquation} X_{t} = x_0 + \int_0^t
F (X_{s}) \,ds + \int_0^t G
(X_{s}) \,dW_{s},
\end{equation}
with coefficients given by
%
\begin{eqnarray}\label{EqDefinitionAverageF}
F(x) &=& \int f(x,y) \mu^x(dy),\qquad G(x) = \sqrt{H(x)},
\nonumber
\\[-8pt]
\\[-8pt]
\nonumber
 H(x)& = &\int h(x,y)
\mu^x(dy),
\end{eqnarray}
where $h(x,y) = gg^*(x,y)$, and $G(x)$ could be in principle any matrix
with square given by $H$, but we choose it to represent the Cholesky
decomposition of the positive semi-definite matrix $H$. It follows that
under appropriate assumptions $X^{\varepsilon} \mathop{\stackrel
{\mathcal{L}}{\longrightarrow}}X$ as $\varepsilon
\rightarrow0$; cf. \citet{pardouxpoisson2003}. The idea behind this
kind of singular perturbation method is that when the difference
between scales is large enough, the dynamics of the system behave as if
the slow scale would be frozen and the ergodic limit of the fast
diffusion would be attained.

However, except for a few particular examples, it is not an easy task
to find explicit expressions for the averages \eqref
{EqDefinitionAverageF}. Naturally, this leads to the question of
designing numerical methods of approximation of the effective equation.
Several methods have been developed for a purely deterministic case;
see, for example, the review \citet{eheterogeneous2007}. Most of them
are based on choosing a macro-solver for the slow scale in which some
information from the fast scale is added via parameters' introduction
to guarantee the correct approximation.

The literature with respect to numerical approximation of the general
stochastic case is, to our knowledge, much more restricted. In \citet
{eanalysis2005} the authors present an algorithm based on the use of
an approximation scheme for the slow scale (e.g., the Euler scheme),
and at each step of the slow scale another scheme is used to solve for
the fast scale contribution; the weak and strong error induced by the
scheme is analyzed when considering the particular case of an ODE with
random coefficients slow scale equation and a stochastic ergodic fast
scale variable [i.e., when $g(x,y)=0$ in \eqref{EqTheSystem}].

In our work we use a similar approach. We focus on approaching
numerically equation \eqref{EqEffectiveEquation}. With this objective
in mind, we propose a \emph{Multi-scale Decreasing Step (MsDS)}
algorithm defined as a composition of an Euler scheme for the slow
scale, the decreasing Euler step algorithm and estimator proposed in
\citet{lambertonrecursive2002} for the ergodic average approximation
at each step, and a Cholesky decomposition for finding the volatility
coefficient.

In order to control the total error approximation of this proposed
algorithm we need to take into account four effects. First, we need an
estimate on the ergodic average approximation at each step. We show
that this control is based on the existence, regularity and control of
the solution of the Poisson equation associated to the fast scale diffusion
%
\begin{equation}
\label{EqPoissonEquation} \mathcal{L}^x_y \phi_\psi(x,y)
= \psi(x,y),
\end{equation}
where
%
\begin{equation}
\label{EqDiffusionOperator} \mathcal{L}^x_y:=\frac{1}{2} \sum
_{i,j = 1}^{d_y} a_{ij}(x,y)
\frac
{\partial^2}{\partial y_i\, \partial y_j} + \sum_{i=1}^{d_y}
b_i(x,y) \frac{\partial}{\partial y_i}
\end{equation}
with $a:=\sigma\sigma^*$, when considering as sources (i.e., the
right-hand side functions) the coefficients $F$ and $H$ centered with
respect to their respective invariant measures.
Second, we need to control the error obtained after performing a
Cholesky decomposition. Then, we have to account for discretization
errors. Finally, we need to control the error propagation which will be
possible under some growth control on the coefficients of the effective
equation.

The MsDS algorithm strongly converges to the exact solution and proves
to be more efficient than a simple Euler scheme for highly oscillating
problems. Moreover, it features a nonstandard C.L.T. property in the
sense that the normalized error distribution converges toward the
solution of an SDE. The coefficients appearing in this normalized error
SDE depend on the solution of the previously mentioned Poisson problem
and are, in general, unknown. Nevertheless, the available explicit
expression for them is valuable for the estimation of confidence
intervals and eases the task of parameter tuning for actual
implementation of the algorithm.

We study as well an \emph{extrapolated MsDS (EMsDS)} version of the
algorithm, differing from the original one in that it uses a
Richardson--Romberg extrapolation of the decreasing step estimator
(i.e., a well-chosen linear combination of the decreasing step Euler
estimator with appropriate parameters) to approach the ergodic
averages. As the MsDS, the EMsDS also features a nonstandard C.L.T.
property and shares the same rate of convergence. However, the
extrapolated version has lower asymptotic complexity and hence higher
asymptotic efficiency than the original one.

\subsection{Outline of the paper}

The organization of the paper is as follows: in Section~\ref
{SecAlgorithm}, we describe the algorithm and state the standing
hypothesis and our main results (strong convergence, limit
distribution). The proof of the main theorem is presented in
Section~\ref{SecMsDSalgorithm} after having reminded some regularity
properties of the effective equation and available results on the
decreasing Euler estimation algorithm in Section~\ref
{SecPreliminaries}. We extend the main results to an extrapolated
version of the algorithm that we introduce and study in Section~\ref
{SecEMSDSalgorithm}. Finally, we perform some numerical studies in
Section~\ref{Secnumresults}. The paper ends with an \hyperref[app]{Appendix}
containing the proof of a couple of technical results.

\section{The MsDS algorithm}
\label{SecAlgorithm}

Let $(\Omega,\F,\prob)$ be a probability space and $W$ be an
$\F$-adapted\vspace*{1pt} Brownian motion. Suppose we are given an independent
probability space $(\tilde{\Omega},\tilde{\F},\tilde{\prob})$ and a
family of independent Brownian motions $\tilde{W}^{q}, q\in\mathds{Q}$
with an associated filtration $\tilde{\F}^{q}_t:=\sigma\{ \tilde
{W}^{q}_s, s\leq t \}$. Define the extended space $(\bar{\Omega},\bar
{\F
},\bar{\F}_t, \bar{\prob})$ by
\begin{eqnarray*}
\bar{\Omega}&:= &\Omega\times\tilde{\Omega},\qquad \bar{\prob}(d\omega ,
d\bar{\omega}) =
\prob(d\omega) \tilde{\prob}(d\bar{\omega}),
\\
\bar{\F}&:=&\F\otimes\tilde{\F},\qquad \tilde{\F}^q_t:= \bigvee
_{q\in
\mathds{Q}; q \leq t} \tilde{F}^q_\infty,\qquad
\bar{F}_t:=\F_t \vee\tilde{\F}^q_t.
\end{eqnarray*}

Such extended space will be useful for treating independently the noise
coming from the Brownian in the effective diffusion and the one related
to the approximation of the ergodic diffusion averages. Consider the
decreasing step Euler algorithm introduced in \citet
{lambertonrecursive2002} to approach the invariant measure of a
recursive diffusion.
Let $\{\gamma_k\}_{k\in\NN}$ be a decreasing sequence of steps satisfying:

\renewcommand{\theHypSteps}{${(\mathcal{H}_{\gamma})}$}
\begin{HypSteps}[(On the sequence of steps for the average
estimation algorithm)]\label{HypSteps}
\begin{longlist}[(iii)]
\item[(i)] $\gamma_k>0$ for all $k$;
\item[(ii)] $\gamma_k$ is a sequence of decreasing steps with $\lim_{n\rightarrow\infty} \gamma_k = 0$;
\label{EnumHipStepsDecreasing}
\item[(iii)] $\lim_{k\rightarrow\infty}\Gamma_k = \infty; \mbox{ where
}\Gamma_k:=\sum_{j=0}^k \gamma_j$;
\label{EnumHipStepsDivergentSum}
\item[(iv)] $\sum_{k=1}^{\infty} ( \frac{\gamma_k^2}{\Gamma_k}
) < + \infty$.
\end{longlist}
\end{HypSteps}

For any $q\in\QQ$, let $\sqrt{\gamma_{k+1}} U_{k+1}^q:= \tilde
{W}^q_{\Gamma_{k+1}}-\tilde{W}^q_{\Gamma_{k}}$ so that $U_{k+1}$ is a
standard Gaussian vector. Let $y_0 \in\RR^{d_y}$. We define the \emph
{decreasing step Euler approximation of the ergodic diffusion} by
%
\begin{eqnarray}
\label{EqMicroAlgo} \DecrAvg{Y}_{0}^{x,q} & =& y_0,
\nonumber
\\[-8pt]
\\[-8pt]
\nonumber
\DecrAvg{Y}_{k+1}^{x,q} & =&\DecrAvg{Y}_{k}^{x,y_0,q}
+ \gamma_{k+1} b \bigl(x,\DecrAvg{Y}_{k}^{x,q}
\bigr)+\sqrt{\gamma_{k+1}} \sigma \bigl(x, \DecrAvg{Y}_k^{x,q}
\bigr)U^q_{k+1},
\end{eqnarray}
and the \emph{decreasing step average estimator} by
%
\begin{equation}
\DecrAvg{F}^{k}(x,q) = \frac{1}{\Gamma_k}\sum
_{j=1}^{k} \gamma_j f \bigl(x,
\DecrAvg{Y}^{x,q}_{j-1} \bigr). \label{EqAverageEstimator}
\end{equation}
The idea behind the particular form of estimator \eqref
{EqAverageEstimator} is to take advantage of the ergodicity of the
diffusion: the long-term time average approaches the invariant measure
of the diffusion. Note that the estimator can also be written
recursively as
\[
\DecrAvg{F}^{0}(x,q) = 0;\qquad \DecrAvg{F}^{k}(x,q) =
\DecrAvg{F}^{k-1}(x,q) + \frac{\gamma_{k}}{\Gamma_k} \bigl(f \bigl(x,
\DecrAvg{Y}_{k-1}^{x,q} \bigr)-\DecrAvg{F}^{k-1}(x,q)
\bigr).
\]

Evidently, using the same ergodic average argument, it is also possible
to use a uniform step estimator of the type $k^{-1}\sum_{j=1}^{k}
\gamma
_j f(x,\DecrAvg{Y}^{x,q}_{j-1})$
as studied, for example, in \citet{Talaysecond-order1990}. The main
difference between both estimators appears in the type of error that
they generate. The uniform step estimator induces two types of errors
coming from the truncation of the series and the fact that the \emph
{ergodic limit of the approached sequence is not the ergodic limit of
the original diffusion}. In contrast, the decreasing Euler scheme
estimator eliminates the asymptotic gap between the invariant law of
the continuous equation and that of its discretization; see \citet
{lambertonrecursive2002}. Moreover, the decreasing step method
features a kind of ``error expansion'' [as shown in \citet
{lemaireestimation2005}] when applied to a certain family of
functions. These properties are important to show the limit properties
of our algorithm and to deduce the extrapolated version.

We should remark that we have chosen to work with a simplified version
of the algorithm in \citet{lambertonrecursive2002}: its more general
version allows the use of different sequences for the Euler scheme step
and for the weights in the average.

With this estimator in hand we can define an Euler scheme to approach
our effective diffusion. Assuming a time horizon $T$, for $n\in\NN^*$
we put $t_k = Tk/n$, so that the Euler scheme will be given by
\[
\check{X}_{t_{k+1}}^{n} = \check{X}_{t_{k}}^{n}
+ \DecrAvg{F}^{M(n)} \bigl(\check{X}_{t_k}^{n},t_k
\bigr) \Delta t_{k+1} + \DecrAvg{G}^{M(n)} \bigl(
\check{X}_{t_k}^{n},t_k \bigr) \Delta
W_{k+1},
\]
where $\DecrAvg{F}^M$ is defined in \eqref{EqAverageEstimator} and
$\DecrAvg{G}^{M}(x,q)$ is defined in two steps: First we find
$\DecrAvg
{H}^{M}(x,q)$ using the decreasing step algorithm as in \eqref
{EqAverageEstimator} [recall that $h(x,y)=g^*g (x,y)$], and then we
perform a Cholesky decomposition on it to find  $\DecrAvg
{G}^{M}(x,q)=\sqrt{\DecrAvg{H}^{M}(x,q)}$. Note that the number of
steps in the decreasing Euler estimator, $M$, is expressed as a
function of the number of steps in the Euler scheme for the slow scale
$n$. The form of $M(n)$ will be clear from the main theorems.

It will be easier to work mathematically with a continuous
interpolation of the Euler approximation. Let us denote by $\tbar(n) =
\lfloor nt\rfloor/n$. We will usually omit the explicit dependence on
$n$ and write $\tbar$ when clear from the context. The continuous Euler
approximation is then given by
%
\begin{equation}
\label{Eqcontinuousapproximated} \ApproxVar{X}^{n}_t = x_0 +
\int_0^t \DecrAvg{F}^{M(n)} \bigl(
\ApproxVar{X}^{n}_{\tbar[s]},\tbar[s] \bigr) \,ds + \int
_0^t \DecrAvg{G}^{M(n)} \bigl(
\ApproxVar{X}^{n}_{\tbar[s]},\tbar[s] \bigr) \,dW_s,
\end{equation}
that is, a linear interpolation from the discrete Euler scheme.
Clearly, at times $t_k$ the continuous Euler coincide with the Euler
algorithm. All our results will be derived for the continuous version
of the algorithm.

\subsection{Standing hypothesis and main result} Let us introduce the
assumptions under which our main results follow.

\renewcommand{\theHypSlowScale}{${(\mathcal{H}_{s.s.})}$}
\begin{HypSlowScale}[(On the slow-scale coefficients)]\label{HypSlowScale}
\begin{longlist}[(iii)]
%
\item[(i)] Lipschitz in x: There exist constants $K, m$ such that for
all $x,x' \in\RR^{d_x}$ and $y \in\RR^{d_y}$,
\[
\bigl|f(x,y)-f \bigl(x',y \bigr)\bigr|+\bigl|g(x,y)-g \bigl(x',y
\bigr)\bigr|\leq K |y|^{m}\bigl|x-x'\bigr|;
\]
\item[(ii)] regularity: $f, h$ belong to $C_{b,p}^{2,r^y}$ for some
$r^y>3$, where the subindex $b,p$ means the derivatives $\partial
_x^i\partial_y^j$ for $0\leq i\leq2$ and $0\leq j \leq r^y-i$ are
bounded in $x$ and polynomially bounded in $y$;
\item[(iii)] degeneracy: either $h$ is identically zero, or it is
uniformly nondegenerate, that is, there exists $\lambda^\prime_{-}
\in
\RR^+_*$ such that $\lambda^\prime_- I \leq h(x,y)$.
\end{longlist}
\end{HypSlowScale}

Before giving the standing hypothesis on the fast scale equation,
recall that we have defined the matrix $a(x,y)=\sigma\sigma^*(x,y)$.

\renewcommand{\theHypFastScale}{${(\mathcal{H}_{f.s.})}$}
\begin{HypFastScale}[(On the fast-scale coefficients)]\label{HypFastScale}
\begin{longlist}[(iii)]
%
\item[(i)] $a,b \in C_{b,l}^{2,0}$, that is, they are continuous and
linearly bounded in $y$ and $C^2$ and bounded in $x$.
\item[(ii)] The matrix $a$ is uniformly continuous and uniformly
nondegenerate and bounded, that is, there exist $\lambda_-, \lambda_+
\in\RR^+_*$ such that
\[
\lambda_ - I \leq a(x,y) \leq\lambda_+ I;
\]
\item[(iii)] $\sup_x b(x,y)\cdot y \leq-c_1|y|^2+c_2$, for some
$c_1\in\RR^*_+,c_2\in\RR$.
\end{longlist}
\end{HypFastScale}

The regularity and growth hypothesis contained in \ref{HypSlowScale} are assumed to
control the error propagation. The main goal of imposing conditions on
the fast scale diffusion is to guarantee the existence of an invariant
limit for any possible fixed value of $x$ and a uniform control on its
averages. For this reasons they are quite restrictive: note that
\ref{HypFastScale}(i) implies $\sup_x |b(x,y)| = O(|y|)$
and \ref{HypFastScale}(iii)
deduces $\lim_{|y|\rightarrow\infty}\sup_x b(x,y)\cdot y = -\infty
$, meaning that
the drift has at most linear growth in $y$ and that it is mean
reverting uniformly in $x$. In turn, the ellipticity and nondegeneracy
assumption \ref{HypFastScale}(ii) is helpful to
deduce the uniqueness of the invariant measure.

We are ready to state our main Theorem on the MsDS algorithm. Its proof
is found in Section~\ref{SecMsDSalgorithm}.

\begin{MyTheorem}
\label{TheoShortMainDecreasingStep}
Let $0<\theta<1$, $\gamma_0 \in\RR^+$ and $\gamma_k = \gamma_0
k^{-\theta}$. Let $M_1$ be a positive constant. Assume
\ref{HypFastScale} and
\ref{HypSlowScale}. Define $M(n)$ by
\[
M(n) = \bigl\lceil M_1 n^{{1}/{(1-\theta)}} \bigr\rceil,
\]
then:
\begin{enumerate}[(ii)]
\item[(i)] ODE with random coefficients case [$g(x,y)\equiv0$]:
\begin{enumerate}[(a)]
\item[(a)] (Strong convergence). There exists a constant $K$ such that
\[
{{\mathbb{E}}} \Bigl[{\sup_{0\leq t \leq T}\bigl|X_t -
\ApproxVar{X}_t^{n}\bigr|^2 }\Bigr] \leq K
n^{-2[(1-\theta)\wedge\theta]/(1-\theta)}.
\]
\item[(b)] 
(Limit distribution of the
error). Assume in addition that $r^y \geq7$ and $ \theta\geq1/2$.
Then
\[
n \bigl(X - \ApproxVar{X}^{n} \bigr) =: \zeta^{n}
\Rightarrow\zeta^\infty,
\]
where $\Rightarrow$ denotes convergence in law, and $\zeta^\infty$ is
the solution of an SDE stated explicitly on Theorem~\ref
{TheoAbstractLimitDistribution}.
\end{enumerate}
\item[(ii)] Full SDE case:
\begin{enumerate}[(a)]
\item[(a)] (Strong convergence). There exists a constant $K$ such that
\[
{{\mathbb{E}}} \Bigl[{\sup_{0\leq t \leq T}\bigl|X_t -
\ApproxVar{X}_t^{n}\bigr|^2 }\Bigr] \leq K
n^{-[(1-\theta)\wedge2\theta]/(1-\theta)}.
\]
\item[(b)] 
(Limit distribution of the
error). Assume in addition that $r^y \geq7$ and $\theta\geq1/3$.
Then
\[
n^{1/2} \bigl(X - \ApproxVar{X}^{n} \bigr) =:
\zeta^{n} \Rightarrow\zeta^{\infty},
\]
where $\zeta^{\infty}$ is the solution of an SDE stated explicitly on
Theorem~\ref{TheoAbstractLimitDistribution}.
\end{enumerate}
\end{enumerate}
\end{MyTheorem}

Note that we study the mean square error of our approximation algorithm
toward the effective equation. We perform this strong error analysis to
guarantee that the algorithm will be used for applications demanding to
approach functions that depend on the whole trajectory (as in finance).
As will be clear from Theorem~\ref{TheoAbstractLimitDistribution}, the
SDE defining the limit results both for the fully stochastic and the
ODE with random coefficients case are explicitly given in terms of the
invariant law of the ergodic diffusion and are consequently unknown.
Nevertheless, the key point is that, being explicit, they might be
estimated numerically for practical purposes.

We have announced an extrapolated version of the algorithm. Given that
its proper introduction requires a further understanding of the basic
algorithm, we postpone the presentation to Section \ref{SecEMSDSalgorithm}.

\section{Preliminaries}
\label{SecPreliminaries}

In this section we present the main tools needed to analyze the
presented algorithm.

Let us start by stating properly the stochastic approximation theorem
we mentioned in the \hyperref[sec1]{Introduction} and that justifies
the relation
between the effective equation \eqref{EqEffectiveEquation} and the
original strongly oscillating system \eqref{EqTheSystem}.

\begin{MyTheorem}[{[Theorem~4 in \citet{pardouxpoisson2003}]}]
\label{IntroTheoPardouxLimit}
Let $b,\sigma,  f,g$ be defined as in \eqref{EqTheSystem} and $a=\sigma
\sigma^*$. Assume we have a recurrence condition of the type $\lim_{|y|\rightarrow\infty} b(x,y)\cdot y=-\infty$, and that the matrix
``$a$'' is nondegenerate and uniformly elliptic. Assume that $a,b \in
C_b^{2,1+\alpha}$, and that $f, g$ are Lipschitz with respect to the
$x$ variable uniformly in $y$ and have at most polynomial growth in $y$
and linear growth in $x$.

Then, for any $T>0$, the family of processes $\{X_t^\varepsilon, 0\leq t
\leq T\}_{0< \varepsilon\leq1}$ is weakly relatively compact in
$C([0,T]; \RR^l)$. Any accumulation point $X$ is a solution of the
martingale problem associated with the operator $\bar{\mathcal{L}}$.

If moreover, the martingale problem is well posed, then $X^\varepsilon
\mathop{\stackrel{\mathcal{L}}{\longrightarrow}}X$, where $X$ is
the unique (in law) diffusion process with
generator $\bar{\mathcal{L}}$.
\end{MyTheorem}

It is worth mentioning that the actual framework of Pardoux and
Vertennikov's statement includes the case in which there is an
$\varepsilon
^{-1}$ order term in the slow variable, which complicates the proof
with respect to the framework we present here. Note that under the
standing hypothesis, the martingale problem is well posed and $X$ in
the theorem is the unique solution to \eqref{EqEffectiveEquation}.

\subsection{A priori estimates}

An important result is related to some a priori estimates valid for
general SDEs. Since they are quite standard, we will state the result
without giving the details of the proof.

\begin{MyProposition}
\label{PropAPriori}
Let
%
\begin{equation}
\vartheta_{t} = \vartheta_0 + \int_0^t
V_1 ( \vartheta_s,s ) \,ds + \int_0^t
V_2 ( \vartheta_s,s) \,dW_{s}, \label{EqGeneralSDEApriori}
\end{equation}
where $V_1, V_2$ are adapted random functions.
\begin{enumerate}[(ii)]
\item[(i)] 
For all $\alpha\geq2$,
\begin{eqnarray*}
&&{{\mathbb{E}}} \Bigl[{\sup_{ 0 \leq t \leq T} |\vartheta_t|^\alpha}
\Bigr] \\
&&\qquad \leq K_\alpha{{\mathbb{E}}} \bigl[{|\vartheta_0|^\alpha}
\bigr] + K( \alpha,T) \int_0^T \bigl( {{
\mathbb{E}}} \bigl[{\bigl| V_1( \vartheta_s,s )
\bigr|^\alpha}\bigr] + {{\mathbb{E}}} \bigl[{\bigl| V_2(
\vartheta_s,s ) \bigr|^\alpha}\bigr] \bigr) \,ds
\\
&& \qquad \leq K_\alpha{{\mathbb{E}}} \bigl[{|\vartheta_0|^\alpha}
\bigr]\\
&&\qquad\quad{} + K'(\alpha,T) \Bigl( \sup_{
0 \leq t \leq T} {{
\mathbb{E}}} \bigl[{\bigl| V_1 (\vartheta_t,t)
\bigr|^\alpha}\bigr] + \sup_{ 0
\leq
t \leq T} {{\mathbb{E}}} \bigl[{\bigl|
V_2 (\vartheta_t, t) \bigr|^\alpha}\bigr] \Bigr).
\end{eqnarray*}
\item[(ii)] 
Assume that $\forall
\alpha\geq2$,
\[
{{\mathbb{E}}} \bigl[{\bigl|V_1(\vartheta_t,t)\bigr|^\alpha}
\bigr] + {{\mathbb{E}}} \bigl[{\bigl|V_2(\vartheta_t,t)\bigr|^\alpha}
\bigr] \leq K \bigl(1+ {{\mathbb{E}}} \bigl[{|\vartheta_t|}\bigr]^{\alpha}
\bigr).
\]
Then:
\begin{enumerate}[(a)]
\item[(a)] 
for $t\in[0,T]$ and
$\alpha\geq2$, ${{\mathbb{E}}} [{|\vartheta_t|^{\alpha}}]\leq
K(\alpha,T)$;
\item[(b)] 
for $\alpha\geq
2$, ${{\mathbb{E}}} [{\sup_{0\leq s\leq t} |\vartheta_s|^{\alpha}}]
\leq K(\alpha
,T) $
$ { {{\mathbb{P}}}( {\sup_{0 \leq s \leq t}\tau_r \leq t } )}\leq
\frac{K'(\alpha
,t)}{r^{\alpha}}$.
\end{enumerate}
\end{enumerate}
\end{MyProposition}

\subsection{Cholesky decomposition}

The Cholesky decomposition of a positive definite matrix consists of
expressing this matrix as the product of a lower triangular matrix and
its conjugate transpose. A stability analysis of this procedure is a
key point in our analysis for the SDE case behavior of our algorithm.

Recall that we denote by $|\cdot|$ the induced operator norm. Let us
denote by $\|\cdot\|_F$ the Frobenius norm. Recall that if $H$ is a
$d\times d$ matrix,
%
\begin{equation}
\label{EqNormEquivalence} |H| \leq\|H\|_F \leq\sqrt{d}|H|.
\end{equation}

\begin{MyTheorem}[{[Theorem~1.1 in \citet{sunperturbation1991}]}]
\label{TheoControlCholesky}
Let $H$ be a $d \times d$ positive definite matrix with Cholesky
factorization $H=GG^*$. If $\Delta H$ is a $d\times d$ symmetrical
matrix satisfying $ |H^{-1}| \|\Delta H\|_F < 1/2$, then there is a
unique Cholesky factorization $H+\Delta H = (G+\Delta G)(G+\Delta G)^*$ and
%
\begin{equation}
\frac{\|\Delta G\|_F}{|G|} \leq\sqrt{2} \frac{\kappa\kappa
_2(H)}{1+\sqrt{1-2\kappa_2(H)\kappa}},
\end{equation}
where $\kappa= |\Delta H\|_F|H|^{-1}$ and $\kappa_2(H)=|H||H^{-1}|$.
\end{MyTheorem}

Theorem~\ref{TheoControlCholesky} gives a control on the sensitivity
of the Cholesky procedure. In Lemma~\ref{LemErrorCholesky} we study
the propagation effect at each stage of the Cholesky factorization to
say a little bit more on the particular form of the error. Its proof is
given in Appendix~\ref{SecAnnexCholesky}.

\begin{MyLemma}
\label{LemErrorCholesky}
Suppose the hypothesis of Theorem~\ref{TheoControlCholesky} holds. Then
\begin{eqnarray*}
\Delta G_{i,i} &=& \frac{\Delta H_{i,i} - 2 \sum_{k=1}^{i-1}\Delta
G_{i,k} G_{i,k}}{2G_{i,i}} + O \bigl(|\Delta
H|^2 \bigr),
\\
\Delta G_{i,j} &=& \frac{\Delta H_{i,j} -G_{i,j}\Delta G_{j,j}
\sum_{k=1}^{j-1} (\Delta G_{j,k} G_{i,k}+\Delta G_{i,k} G_{j,k})}{G_{j,j}} + O \bigl(|\Delta H|^2
\bigr)
\end{eqnarray*}
for $i>j$.
\end{MyLemma}

Lemma~\ref{LemErrorCholesky} gives a first order approximation of the
error matrix $\Delta G$ knowing the perturbation matrix $\Delta H$.
From this lemma, we can deduce on the regularity of the Cholesky
approximation. The following corollary follows from the definition of
$H$ and Lemma~\ref{LemErrorCholesky}.

\begin{MyCorollary}
\label{CorRegularityG}
Let $H\dvtx\RR^{d} \rightarrow M^{d \times d}$ be $C_b^2$ and
nondegenerate [in the sense given in Hypothesis
\ref{HypSlowScale}]. Then $G$ is
also $C_b^2$ and nondegenerate.
\end{MyCorollary}

\subsection{Decreasing step Euler algorithm} In this section we present
some control and error expansion results valid for the decreasing step
Euler algorithm. The results here presented are found in \citet
{lambertonrecursive2002} or in the Ph.D. thesis of \citet
{lemaireestimation2005}.

A first interesting property is that the sequence of estimators defined
in \eqref{EqAverageEstimator} converges almost surely to the ergodic
average for any fixed $x$.

\begin{MyProposition}
\label{PropConvergenceEstimator}
Assume \ref{HypFastScale}, and let $\psi\dvtx\RR^{d_x}\times
\RR^{d_y}\rightarrow
\RR$,
and suppose that $\psi(x,y) \leq C(x)(1+|y|^\pi)$. Let $\DecrAvg
{\Psi
}^M(x,q)$ be defined as in \eqref{EqAverageEstimator}.
Then, for any $x\in\RR^{d_x}, q\in\QQ$,
\[
\DecrAvg{\Psi}^M(x,q) \tendsas\int\psi(x,y)\mu^x(dy)\qquad
\mbox{as } M\rightarrow\infty,
\]
where $\mu^x$ is the invariant measure of \eqref{EqErgodicDiffusion}.
\end{MyProposition}

\begin{pf}
\ref{HypFastScale} imply that $V(y):= 1 + |y|^2$ is a uniformly
in $x$ function
satisfying the hypothesis of Theorem~1 in \citet
{lambertonrecursive2002}, from which the claim follows.
\end{pf}

We have as well a control on the moments of any order of $\DecrAvg{Y}^{x,q}_k$.

\begin{MyProposition}
\label{PropApproxFiniteMomentstime}
Let $\pi>0$ and let $\DecrAvg{Y}_k^{x,q}$ be given by \eqref
{EqMicroAlgo}. Then there exists a constant $K_\pi$ given only by $\pi
$, $\lambda_-$, $\lambda_+$ and $\gamma_0$ such that for all $x\in
\RR
^{d_x}$ and $q\in\QQ$,
\[
\sup_{i \in\NN} {{\mathbb{E}}} \bigl[{\bigl|\DecrAvg{Y}_i^{x,q}\bigr|^\pi}
\bigr] < K_\pi.
\]
Moreover, for every $\pi> 1$,
\[
\sup_{M\in\NN} \Biggl(\frac{1}{\Gamma_M}\sum
_{i=1}^{M} \gamma_i \bigl|
\DecrAvg{Y}_i^{x,q}\bigr|^\pi \Biggr) < + \infty.
\]
\end{MyProposition}

\begin{pf}
By Lemma~2 in \citet{lambertonrecursive2002} given that $U_k^q$ has
moments of any order and $V(y) = |y|^2 +1$ satisfies the needed
hypothesis uniformly in $x$, we get that for any $\pi\geq1$ and $q\in
\QQ$,
\[
\sup_{i \in\NN} {{\mathbb{E}}} \bigl[{\bigl|\DecrAvg{Y}_i^{x,q}\bigr|^{2\pi}}
\bigr] \leq\sup_{i
\in\NN
} {{\mathbb{E}}} \bigl[{V \bigl(
\DecrAvg{Y}_i^{x,q} \bigr)^{\pi}}\bigr] <
K_{\pi}.
\]
The extension to all $\pi>0$ is straightforward.

The second claim follows from Theorem~3 in \citet{lambertonrecursive2002}.
\end{pf}

Proposition~\ref{PropErrorExpansion} is an adaptation of a result
appearing in the Ph.D. thesis \citet{lemaireestimation2005}. The proof
comes from performing a Taylor expansion and reordering the terms in a
proper way. For the statement, we introduce in addition to the sequence
$\{\gamma_k\}_{\{k\in\NN^*\}}$ a new sequence that we denote by $\{
\eta
_k\}_{\{k\in\NN^*\}}$ (that may be taken equal to the former). This
added flexibility will be useful in the following, in particular to
prove Proposition~\ref{PropControlDecreasingDifferentStep}. We may
interpret Proposition~\ref{PropErrorExpansion} as an error expansion
result. Indeed if we fix $\eta_k=\gamma_k$ satisfying
\ref{HypSteps}, then
we will have an explicit expression for the approximation error of the
decreasing Euler algorithm.

\begin{MyProposition}
\label{PropErrorExpansion}
Let $\psi\dvtx\RR^{d_x}\times\RR^{d_y}\rightarrow\RR$. Under the
assumptions of Proposition~\ref{PropConvergenceEstimator}, suppose
that for each $x\in\RR^{d_x}$ there exists $\phi_\psi^x\dvtx\RR
^{d_y}\rightarrow\RR$ solution of the centered Poisson equation
%
\begin{equation}
\label{EqCenteredPoisson} \mathcal{L}_y^x \phi^x_\psi(y)
= \psi(x,y)-\int\psi(x,z)\mu^x(dz).
\end{equation}
Suppose as well for $r\in\NN$, $r\geq2$, that $\phi^x_\psi$ is
$C^{r}$ in the $y$-variable uniformly in~$x$, and $D^r \phi_\psi$ is
Lipschitz in $y$ uniformly in $x$. Let $\gamma_k$ and $\eta_k$ be two
decreasing sequences with\vspace*{1pt} $\gamma_k \rightarrow0$, $\eta
_k\rightarrow
0$, $\Gamma_k=\sum_{1\leq j\leq k}\gamma_k$, $\mathcal{H}_k=\sum_{1\leq
j\leq k}\eta_k$. Let $\DecrAvg{Y}_k^{x,q}$ be defined as in \eqref
{EqMicroAlgo} (with step sequence $\gamma_k$). Then
\[
\sum_{k=1}^{M} \eta_k \biggl(
\psi \bigl(x,\DecrAvg{Y}^{x,q}_{k-1} \bigr) -\int\psi(x,z)
\mu^x(dz) \biggr) = A^{0}_{\psi,M}-N_{\psi,M}
- \sum_{i=2}^{r} A^{i}_{\psi,M}
-Z^{r}_{\psi,M},
\]
where
%
\begin{eqnarray}
A^{0}_{\psi,M}  (x,q)&:=& \sum_{k=1}^{M}
\frac{\eta_k}{\gamma_k} \bigl[ \phi_\psi^x \bigl(
\DecrAvg{Y}^{x,q}_{k} \bigr)-\phi_\psi^x
\bigl(\DecrAvg{Y}^{x,q}_{k-1} \bigr) \bigr], \label{EqA0M}
\\
N_{\psi,M}  (x,q)&:=&\sum_{k=1}^{M}
\frac{\eta_k}{\sqrt{\gamma_k}} \bigl\langle D_y \phi_\psi^x
\bigl(\DecrAvg{Y}^{x,q}_{k-1} \bigr), \sigma \bigl(x,\DecrAvg
{Y}^{x,q}_{k-1} \bigr)U_k^q \bigr
\rangle,\label{EqNM}
\\
A^{2}_{\psi,M}  (x,q)&:=& \frac{1}{2}\sum
_{k=1}^{M}\eta_k \bigl[ D^2
\phi_\psi^x \bigl(\DecrAvg{Y}^{x,q}_{k-1}
\bigr)\cdot \bigl(\sigma \bigl(x,\DecrAvg{Y}^{x,q}_{k-1}
\bigr)U^q_k \bigr)^{\otimes2}\label{EqA2M}
\nonumber
\\[-8pt]
\\[-8pt]
\nonumber
& &\hspace*{39pt}{} - \operatorname{Tr}( D^2 \phi_\psi^x
\bigl(\DecrAvg{Y}^{x,q}_{k-1} \bigr) \bigl(\sigma^*\sigma
\bigl(x,\DecrAvg{Y}^{x,q}_{k-1} \bigr) \bigr) \bigr],
\\
A^{i}_{\psi,M} (x,q)&:=& \sum_{k=1}^{M}
\eta_k\gamma_k^{i/2-1} v_\psi
^{i,r} \bigl(x,\DecrAvg{Y}^{x,q}_{k-1},U_k^q
\bigr) \label{EqAiM}
\end{eqnarray}
for $i = 3,\dots, r$ with
\[
v_\psi^{i,r}(x,y,z) = \sum_{j\geq i/2}^{i\wedge r}
\pmatrix{j
\cr
i-j}\frac
{1}{j!} D_y^j
\phi_\psi^x(y) \cdot \bigl\langle b(x,y)^{\otimes(i-j)},
\bigl(\sigma(x,y) z \bigr)^{\otimes(2j-i)} \bigr\rangle
\]
and
%
\begin{equation}
\bigl|Z_{\psi,M}^{r}\bigr| (x,q)\leq K \sum
_{k=1}^{M} \eta_k\gamma_k^{{(r-1)}/{2}}
\bigl(1+\bigl|\DecrAvg{Y}^{x,q}_{k-1}\bigr|^{r+1} \bigr) \bigl(1
+ \bigl|U_k^q\bigr| \bigr)^{r+1}.\label{EqZ2M}
\end{equation}
\end{MyProposition}

The average of each expansion term will play an important role in our
analysis, so that we will present a special notation for them. Indeed, let
\begin{eqnarray}\label{EqDefinitionvbar}
&&\bar{v}_\psi^{i,r} (x,y)\nonumber \\
&&\qquad:= {{\mathbb{E}}}
\bigl[{v_\psi^{i,r} \bigl(x,y,U_1^0
\bigr)}\bigr]
\\
&&\qquad = \sum_{j\geq i/2}^{i\wedge r} \pmatrix{j
\cr
i-j}
\frac
{1}{j!}D^j_y\phi_\psi^x(y)
{{\mathbb{E}}}\bigl[ { \bigl\langle b(x,y)^{\otimes(i-j)}, \bigl(\sigma (x,y)
U_k^q \bigr)^{\otimes(2j-i)} \bigr\rangle} | {\tilde
{F}_{\Gamma_{k-1}}} \bigr]
\nonumber
.
\end{eqnarray}

\begin{MyRemark} \label{RmkZeroExpectation} Consider $A_{\psi
,M}^{2i+1}$ for $i\leq\lfloor(r-1)/2 \rfloor$.
As $2j-2i-1$ is odd for any $j$ integer and given the fact that the odd
powers of a centered Gaussian are centered, we deduce
$\bar{v}_\psi^{2i+1,r} = 0.$ Of course this property transfers to
$A_{\psi,M}^{2i+1}$ so that ${{\mathbb{E}}} [{A_{\psi
,M}^{2i+1}}]=0$, implying in
turn that the terms with an odd index are centered.
\end{MyRemark}

Under some additional hypotheses, Proposition~\ref{PropErrorExpansion}
may be used to obtain an $L_2$ control on the error of the
approximation. For the sake of the presentation, let us denote from now
on
%
\begin{equation}
\label{EqDefGammaAlpha}\Gamma_M^{[r]} = \sum
_{k=1}^M (\gamma_k)^r.
\end{equation}
Note we have in particular $\Gamma_M^{[1]}=\Gamma_M$.

\begin{MyProposition}
\label{PropControlDecreasingDifferentStep}
Under the assumptions of Proposition~\ref{PropErrorExpansion}, let
$\alpha\geq1$. Assume $\{\gamma_k\}$ satisfies
\ref{HypSteps}, and that
$\Gamma_M^{[\alpha]}\rightarrow\infty$, for $\Gamma_M^{[\alpha]}$
defined as in \eqref{EqDefGammaAlpha}. Assume as well that the
solution of the centered Poisson equation $\phi_\psi$ is in
$C_{b,p}^{2,r}$ for $r> 3$. Let $\bar{\Psi}:= \int\psi(x,z)\mu^x(dz)
$, then
\[
{{\mathbb{E}}} \Biggl[{\Biggl| \frac{1}{\Gamma_M^{[\alpha]}}\sum_{k=1}^{M}
\gamma^\alpha_k \bigl(\psi \bigl(x,\DecrAvg{Y}^{x,q}_{k-1}
\bigr) -\bar{\Psi}(x) \bigr) \Biggr|^2}\Biggr] \leq K \frac{ 1+ \Gamma^{[2\alpha
-1]}_M +\Gamma^{[2\alpha]}_M +
(\Gamma
^{[\alpha+1]}_M)^2}{(\Gamma^{[\alpha]}_M)^2}.
\]
\end{MyProposition}

\begin{pf}
We recall first some martingale inequalities. Let $\{a_k\}$ be any
sequence of random tensors. By Cauchy--Schwarz inequality we have that
%
\begin{equation}
\label{EqEstimateSquareGeneral} {{\mathbb{E}}} \Biggl[{\Biggl|\sum_{k=1}^{M}
\gamma_k^p a_{k}\Biggr|^2}\Biggr]
\leq{{\mathbb{E}}} \Biggl[{\Gamma^{[p]}_M \sum
_{k=1}^{M} \gamma_k^p
|a_{k}|^2}\Biggr] = \Gamma^{[p]}_M
\sum_{k=1}^{M} \gamma_k^p
{{\mathbb{E}}} \bigl[{|a_{k}|^2}\bigr].
\end{equation}
Let $\{b_k\}$ be also a sequence of tensors. If $s_0 < s_1 < \cdots<
s_k < \cdots,$ the $\{a_k\}, \{b_k\}$ are $\DecrAvg{\F}_{s_k}^q$
adapted, and for all $k$, ${{\mathbb{E}}}[ {a_k} | {\DecrAvg{\F
}^q_{s_k}} ]={{\mathbb{E}}}[ {b_k} | {\DecrAvg{\F}^q_{s_k}} ]=0$, we
have by martingale properties that
%
\begin{equation}
{{\mathbb{E}}} \Biggl[{ \Biggl\langle\sum_{k=1}^{M}
\gamma_k^p a_{k}, \sum
_{k=1}^{M}\gamma_k^p
b_{k} \Biggr\rangle}\Biggr] = \sum_{k=1}^{M}
\gamma_k^{2p}{{\mathbb{E}}} \bigl[{\langle
a_{k}, b_{k} \rangle}\bigr] \label
{EqEstimateProductCentered}
\end{equation}
and in particular,
%
\begin{equation}
{{\mathbb{E}}} \Biggl[{\Biggl|\sum_{k=1}^{M}
\gamma_k^p a_{k}\Biggr|^2}\Biggr] = \sum
_{k=1}^{M} \gamma_k^{2p}{{
\mathbb{E}}} \bigl[{| a_{k} |^2}\bigr]. \label{EqEstimateSquareCentered}
\end{equation}
Now, take the error expansion in Proposition~\ref{PropErrorExpansion}
with $r=3$, and let $\eta_k = \gamma_k^\alpha$. By Abel's
transformation, using convexity, estimate \eqref
{EqEstimateSquareGeneral}, the regularity properties of $\phi_\psi$
and Proposition~\ref{PropApproxFiniteMomentstime}, we get
%
\begin{eqnarray}
\label{EqControlA0} %
&& {{\mathbb{E}}} \bigl[{\bigl|A_{\psi,M}^{0}(x,q)
\bigr|^2}\bigr] \nonumber\\
&&\qquad={{\mathbb{E}}} \Biggl[{\Biggl| \sum
_{k=1}^{M} \gamma_k^{\alpha-1}
\bigl[ \phi_\psi^x \bigl(\DecrAvg{Y}^{x,q}_{k}
\bigr)-\phi_\psi ^x \bigl(\DecrAvg {Y}^{x,q}_{k-1}
\bigr) \bigr]\Biggr|^2}\Biggr]
\nonumber\\
&&\qquad = {{\mathbb{E}}} \Biggl[\Biggl| \gamma_M^{\alpha-1}
\phi_\psi^x \bigl(\DecrAvg{Y}^{x,q}_{M}
\bigr) - \gamma_0^{\alpha-1} \phi_\psi ^x
\bigl(\DecrAvg{Y}^{x,q}_{0} \bigr) \nonumber\\
&&\hspace*{14pt}\qquad\quad{}+ \sum
_{k=1}^{M-1} \bigl[ \bigl( \gamma_k^{\alpha-1}-
\gamma_{k+1}^{\alpha-1} \bigr) \phi _\psi^x
\bigl(\DecrAvg{Y}^{x,q}_{k} \bigr) \bigr]\Biggr|^2
\Biggr]
\\
& &\qquad\leq 3 {{\mathbb{E}}} \bigl[{\bigl|\gamma_M^{\alpha-1}
\phi_\psi^x \bigl(\DecrAvg{Y}^{x,q}_{M}
\bigr) \bigr|^2}\bigr] + 3 {{\mathbb{E}}} \bigl[{\bigl| \gamma
_0^{\alpha-1}\phi_\psi^x \bigl(
\DecrAvg{Y}^{x,q}_{0} \bigr)\bigr |^2}\bigr] \nonumber\\
&&\qquad\quad{}+ 3 {{
\mathbb{E}}} \Biggl[{\Biggl| \sum_{k=1}^{M-1}
\bigl[ \bigl( \gamma _k^{\alpha-1}-\gamma_{k+1}^{\alpha-1}
\bigr) \phi_\psi^x \bigl(\DecrAvg{Y}^{x,q}_{k}
\bigr) \bigr]\Biggr|^2}\Biggr]
\nonumber\\
&&\qquad \leq K \Biggl[ \bigl(\gamma_M^{\alpha-1} \bigr)^2
+1 + \Biggl(\sum_{k=1}^{M-1} \bigl(
\gamma_k^{\alpha-1}-\gamma_{k+1}^{\alpha-1} \bigr)
\Biggr)^2 \Biggr] \leq K.\nonumber %
\end{eqnarray}

Moreover, using the fact that the terms are centered from Remark~\ref
{RmkZeroExpectation}, equation~\eqref{EqEstimateSquareCentered} and
the finite moments of the Brownian increments imply
%
\begin{eqnarray}\qquad
{{\mathbb{E}}} \bigl[{\bigl|N_{{\psi},M} (x,q )\bigr|^2}\bigr] &=& \sum
_{k=1}^M \gamma_k^{2\alpha
-1}{{
\mathbb{E}}} \bigl[{\bigl| \bigl\langle\sigma^* D_y \phi_{\psi}
\bigl(x, \DecrAvg{Y}^{x,q}_{k-1} \bigr), U_k^q
\bigr\rangle\bigr|^2}\bigr] \leq K \Gamma_M^{2\alpha-1},\label{EqControlN}
\\
\label{EqControlA2}
{{\mathbb{E}}} \bigl[{\bigl|A_{{\psi},M}^{2}(x,q)\bigr|^2}
\bigr] &\leq& \frac{1}{4} \sum_{k=1}^M
\gamma _k^{2\alpha} {{\mathbb{E}}} \bigl[{\bigl|D^2_y
\phi_\psi \bigl(x,\DecrAvg {Y}^{x,q}_{k-1} \bigr)
\cdot \bigl(\sigma \bigl(x,\DecrAvg {Y}^{x,q}_{k-1}
\bigr)U^q_k \bigr)^{\otimes2} \bigr|^2}\bigr]
\nonumber
\\[-8pt]
\\[-8pt]
\nonumber
&\leq& K\Gamma^{[2\alpha]}_M. %
\end{eqnarray}
More generally, estimate \eqref{EqEstimateSquareCentered} leads to
%
\begin{equation}
\label{EqControlA3} {{\mathbb{E}}} \bigl[{\bigl|A_{{\psi},M}^{3}(x,q)\bigr|^2}
\bigr] = \sum_{k=1}^M
\gamma_k^{2\alpha
+1} {{\mathbb{E}}} \bigl[{\bigl|v_{\psi}^{3,r}
\bigl(x,\DecrAvg {Y}^{x,q}_{k-1},U_k^q
\bigr)\bigr|^2}\bigr] \leq K \Gamma_M^{[2\alpha+1]},
\end{equation}
while by virtue of \eqref{EqEstimateSquareGeneral}, we find as estimate
%
\begin{eqnarray}
\label{EqControlZ3} {{\mathbb{E}}} \bigl[{\bigl|Z_{{\psi},M}^{3}(x,q)\bigr|^2}
\bigr] &\leq& K {{\mathbb{E}}} \Biggl[{\Biggl|\sum_{k=1}^M
\gamma _k^{\alpha+1} \bigl(1+\bigl|\DecrAvg{Y}^{x,q}_{k-1}\bigr|^{4}
\bigr) \bigl(1+\bigl|U_k^q\bigr| \bigr)^4\Biggr|^2}
\Biggr]
\nonumber
\\[-8pt]
\\[-8pt]
\nonumber
&\leq& K \bigl(\Gamma_M^{[\alpha+1]} \bigr)^2.
\end{eqnarray}
On the other hand, from \ref{HypSteps} and given that
$\Gamma_M^{[\alpha
]}\rightarrow\infty$, we have for $M$ large enough that, if $i>j$,
\[
\frac{\Gamma_M^{[i]}}{\Gamma^{[\alpha]}_M}\leq\frac{\Gamma
_M^{[j]}}{\Gamma^{[\alpha]}_M}.
\]
The claim follows from Proposition~\ref{PropErrorExpansion} and \eqref
{EqControlA0}--\eqref{EqControlZ3}.
\end{pf}

\subsection{Ergodic average and Poisson equation} Being basic to our
analysis, we introduce in this section some known properties of the
exact averages and the effective diffusion. These results are studied
in \citeauthor{pardouxpoisson2001}
(\citeyear{pardouxpoisson2001,pardouxpoisson2003}).

Let us start by stating a growth control result proved in \citet
{veretennikovpolynomial1997}.

\begin{MyProposition}\label{PropFiniteMomentsTime} Let $\alpha>0$, and
let $Y_t^{x}$ be the solution of \eqref{EqErgodicDiffusion} with
deterministic initial condition $y_0$ and coefficients satisfying
\ref{HypFastScale}.

Then there exists a constant $K$ given only by $\alpha$, $\lambda_-$,
$\lambda_+$ such that for all $t\geq0$ and $x\in\RR^{d_x}$,
\[
{{\mathbb{E}}} \bigl[{\bigl|Y_t^{x}\bigr|^\alpha}\bigr] <
K \bigl(1+|y_0|^{\alpha+2} \bigr).
\]
\end{MyProposition}

This proposition has a natural corollary.

\begin{MyCorollary}\label{CorFiniteMomentsIM} Under the same
hypothesis of the theorem, for any $\alpha>0$ and all $x\in\RR^{d_x}$,
\[
\int|y|^\alpha\mu^x(dy) < K.
\]
\end{MyCorollary}

\begin{MyLemma}
\label{LemPropertiesAverage}
Let $\psi(x,y)$ be a function satisfying the regularity and growth
conditions in \ref{HypSlowScale}, and let $\Psi(x)=\int\psi
(x,y) \mu^x(dy)$, then
$\Psi(x)$ is~$C^2_b$.
\end{MyLemma}

\begin{pf}
The claim follows from adapting Theorems 3 and 5 in \citet
{veretennikovsobolev2011} to the linear growth case: the needed
equivalent results of convergence in total variation and control of
expectations may be found in \citet{meynstability1993}.
\end{pf}

As it was shown in Proposition~\ref{PropErrorExpansion}, the centered
Poisson equation \eqref{EqCenteredPoisson} plays a special role in
understanding the error expansion of the decreasing Euler algorithm.
Proposition~\ref{PropPropertiesPhi}, which is an adaptation of Theorem~1 in \citet{pardouxpoisson2001} and \citet{veretennikovsobolev2011},
states some sufficient conditions for having the solution of such an
equation when $f$ belongs to a certain family of functions.

\begin{MyProposition}
\label{PropPropertiesPhi}
Consider a function $\psi(x,y)$ satisfying the regularity and growth
conditions in \ref{HypSlowScale}\textup{(i)}, \textup{(ii)} and such that
\[
\int\psi(x,y)\mu^x(dy) = 0 \qquad\forall x.
\]
Assume \ref{HypFastScale}. Then there exists a function $\phi
_\psi(x,y)$,
continuous in $y$ and belonging to the class $\bigcap_{p>1}
W_{p,\mathrm{loc}}^2$ in $y$, such that for every $x\in\RR^{d_x}$:
\begin{enumerate}[(iii)]
\label{EnumSolvesPoisson}
\item[(i)] $\mathcal{L}^x_y \phi_\psi(x,y) = \psi(x,y)$,
\label{EnumIsCentered}
\item[(ii)] $\int\phi_\psi(x,y) \mu^x(dy) = 0 $,
\item[(iii)] $\phi_\psi\in C^{2,r^y}_{b,p}$.
\end{enumerate}
This function is the unique solution up to an additive constant of the
Poisson equation on the class of continuous and $\bigcap_{p>1}
W_{p,\mathrm{loc}}^2$ functions in $y$ which are locally bounded and grow at
most polynomially in $|y|$ as $|y|\rightarrow\infty$. Moreover, it has
the representation
\[
\phi_\psi(x,y) = - \int_0^\infty
\mathbb{E}_{x,y} \bigl( \psi \bigl(x,Y_t^x
\bigr) \bigr) \,dt.
\]
\end{MyProposition}

\section{Convergence results for the MsDS algorithm}
\label{SecMsDSalgorithm}

We focus now on the study of the MsDS algorithm. First, we show that
the proposed approximated coefficients (by means of Decreasing Euler
step and Cholesky procedures) satisfy a growth control and error
control properties. As a consequence, we will conclude on some
regularity property of the approximated diffusion \eqref
{Eqcontinuousapproximated} and show its strong convergence toward
\eqref{EqEffectiveEquation}. Then we will study the limit error
distribution property.

\subsection{Existence, uniqueness, continuity}

From Hypotheses \ref{HypSlowScale}, \ref{HypFastScale},
Proposition~\ref
{PropFiniteMomentsTime} and Proposition~\ref{PropAPriori}, it follows
that there exists a unique solution to equation (\ref
{EqEffectiveEquation}), and that it has a continuous modification. We
show the defined approximation has the same properties.

Proposition~\ref{PropControlDecreasingStep} uses the results of
Section~\ref{SecPreliminaries} to show that, under the standing
hypothesis, the coefficients of the approximated diffusion have finite
moments of any order, and that its error with respect to the exact
coefficients decrease as a power of the number of steps $n$.

\begin{MyProposition}
\label{PropControlDecreasingStep}
Assume \ref{HypSlowScale}, \ref{HypFastScale} and
\ref{HypSteps}. Let $\beta_0> 0$, and define
$M(n)$ implicitly by $\Gamma_{M(n)} = C_0 n^{2\beta_0}$, where $C_0$ is
some constant.
\begin{enumerate}[(ii)]
\item[(i)] 
There exist $\phi
_f$ and $\phi_h$ solutions of the centered Poisson equations:
\begin{itemize}
\item$\mathcal{L}_y^x \phi_f(x,y) = f(x,y) - \int f(x,y')\mu^x(dy')$;
\item$\mathcal{L}_y^x \phi_h(x,y) = h(x,y) - \int h(x,y')\mu^x(dy')$.
\end{itemize}

\item[(ii)] 
Let
%
\begin{equation}
\label{EqDefinitionLbar} \varsigma:= \min_{l\geq4, i=1,\ldots, d} \bigl(
\bar{v}^{l,r^y}_{F^i} \neq0 \bigr)\wedge\min_{l\geq4, i,j=1,\ldots, d}
\bigl( \bar{v}^{l,r^y}_{H^{i,j}} \neq0 \bigr) \wedge
\bigl(r^y+1 \bigr)
\end{equation}
[with the convention that $\min(\varnothing) = \infty$] and $\bar
{v}^{l,r}_{F^i}, \bar{v}^{l,r^y}_{H^{i,j}}$ defined as in \eqref
{EqDefinitionvbar} applied to $F^1,\ldots,F^{d_x}$, $H^{1,1},\ldots,
H^{d_x,d_x} $. Assume the asymptotic expantion
%
\begin{equation}
\label{EqConditionGammaMainDecreasingStep} \frac{\Gamma^{[\varsigma
/2]}_M}{\Gamma_M} = C_1 n^{-\beta_1} + o
\bigl(n^{-\beta_1} \bigr),
\end{equation}
for some $\beta_1>0$, and some constant $C_1$, holds. Let
%
\begin{equation}
\label{EqDefinitionBeta} \beta:=\beta_0\wedge\beta_1.
\end{equation}
Then $\DecrAvg{F}^n$ (and resp., $\DecrAvg{H}^n, \DecrAvg
{G}^n:=\sqrt
{\DecrAvg{H}^n}$) satisfies for any $\alpha\in\RR^+$ and
$k=0,\ldots,n$
\[
\cases{ {{\mathbb{E}}} \bigl[{\bigl|\DecrAvg{F}^n(x,t_k)\bigr|
^\alpha}\bigr] \leq K, \vspace*{2pt}
\cr
{{\mathbb{E}}} \bigl[{\bigl|
\DecrAvg{F}^{n}(x,t_k) - F(x)\bigr|^2}\bigr] \leq K
n^{-2\beta}.}
\]
\end{enumerate}
\end{MyProposition}

\begin{MyRemark}
We should understand $\varsigma$ as marking the first nonzero value in
the error expansion of either $\DecrAvg{F}^n$ or $\DecrAvg{H}^n$. It
depends exclusively on the coefficients of the effective and ergodic
diffusion (in particular it does not depend on $n$).
\end{MyRemark}

\begin{MyRemark}
Proposition~\ref{PropControlDecreasingStep} means that we have a rate
of convergence in norm $L_2$ for the coefficient estimators of order
$O(n^{-\beta})$. Since we \emph{choose} $\beta_0$ by taking $M(n)$ as
needed, the actual limit to $\beta$ comes from $\beta_1$. But of
course, increasing $\beta_0$ implies growing $M$ faster as a function
of $n$, increasing the algorithm's cost.
\end{MyRemark}

\begin{pf*}{Proof of Proposition~\ref{PropControlDecreasingStep}}
Note first that (i) follows
from \ref{HypSlowScale} and Proposition~\ref{PropPropertiesPhi}.

We prove (ii). By Jensen's
inequality and Proposition~\ref{PropApproxFiniteMomentstime}, we have
for every $\alpha\geq1$ and $n$ big enough,
\begin{eqnarray*}
{{\mathbb{E}}} \bigl[{\bigl|\DecrAvg{F}^n(x,q)\bigr|^\alpha}\bigr] &=&
{{\mathbb{E}}} \Biggl[{\Biggl|\frac{1}{\Gamma_M} \sum_{k=1}^M
\gamma_k f \bigl(x,\DecrAvg {Y}^{x,q}_{k-1}
\bigr)\Biggr|^\alpha}\Biggr] \\
&\leq&{{\mathbb{E}}} \Biggl[{\frac
{1}{\Gamma_M}\sum
_{k=1}^M \gamma_k \bigl|f \bigl(x,
\DecrAvg {Y}^{x,q}_{k-1} \bigr)\bigr|^\alpha}\Biggr] \leq K,
\end{eqnarray*}
and similarly for every $\alpha\geq2$,
\[
{{\mathbb{E}}} \bigl[{\bigl|\DecrAvg{G}^n(x,q)\bigr|^\alpha}\bigr] =
{{\mathbb{E}}} \bigl[{\bigl|\DecrAvg {H}^n(x,q)\bigr|^{\alpha/2}}\bigr]
\leq K,
\]
since $|G|^2=|H|$. The result extends trivially to every $\alpha> 0$.

It remains to prove the error control. We obtain an expansion of order
$r^y$ in Proposition~\ref{PropErrorExpansion}. We can bound the first
terms as we did in Proposition~\ref
{PropControlDecreasingDifferentStep} by taking $\gamma_k=\eta_k$ for
all $k=1,\ldots,M$ (i.e., taking $\alpha=1$ in the statement of
Proposition~\ref{PropControlDecreasingDifferentStep}). More generally,
from the definition of $\varsigma$ in \eqref{EqDefinitionLbar}, we
have that for every $l<\varsigma$ or $l$ odd $\bar
{v}^{l,r^y}_{F^{i}}(x,y)=0$, \eqref{EqEstimateSquareCentered} leads to
%
\begin{equation}
\label{EqControlAicentered} {{\mathbb{E}}} \bigl[{\bigl|A_{{F^{i}},M}^{l}(x,q)\bigr|^2}
\bigr] = \sum_{k=1}^M
\gamma_k^l {{\mathbb{E}}} \bigl[{\bigl|v_{F^i}^{l,r^y}
\bigl(x,\DecrAvg {Y}^{x,q}_{k-1},U_k^q
\bigr)\bigr|^2}\bigr] \leq K \Gamma_M^{[l]},
\end{equation}
while for even $l$ with $l\geq\varsigma$, by virtue of \eqref
{EqEstimateSquareGeneral}, we find as estimate
%
\begin{equation}
\label{EqControlAinoncentered} \qquad{{\mathbb{E}}}
 \bigl[{\bigl|A_{{F^{i}},M}^{l}
(x,q)\bigr|^2}
\bigr] \leq \Gamma^{l/2}_M \sum
_{k=1}^M \gamma_k^{l/2} {{
\mathbb{E}}} \bigl[{\bigl|v_{F^i}^{l,r^y} \bigl(x,\DecrAvg
{Y}^{x,q}_{k-1},U_k^q
\bigr)\bigr|^2}\bigr] \leq K \bigl(\Gamma_M^{[l/2]}
\bigr)^2.
\end{equation}
Likewise,
%
\begin{eqnarray}
\label{EqControlZi} %
&&{{\mathbb{E}}} \bigl[{\bigl|Z_{{F^{i}},M}^{r^y}(x,q)\bigr|^2}
\bigr]\nonumber\\
&&\qquad \leq K {{\mathbb{E}}} \Biggl[{\Biggl|\sum_{k=1}^M
\gamma_k^{r^y+1/2} \bigl(1+\bigl|\DecrAvg{Y}^{x,q}_{k-1}\bigr|^{r^y+1}
\bigr) \bigl(1+\bigl|U_k^q\bigr| \bigr)^{r^y+1}\Biggr|^2}
\Biggr]
\\
 &&\qquad \leq K \bigl(\Gamma_M^{[r^y+1/2]} \bigr)^2.\nonumber
\end{eqnarray}
Note that estimates \eqref{EqControlAicentered} and \eqref
{EqControlAinoncentered} are uniform in $x$. On the other hand, from
\ref{HypSteps}, we have for $M$ big enough and $l\leq
r^y$ that
\[
1\geq\frac{\Gamma_M^{[2]}}{\Gamma_M}\geq\frac{\Gamma
_M^{[3]}}{\Gamma
_M}\geq\cdots\geq\frac{\Gamma_M^{[l]}}{\Gamma_M}.
\]
Hence from Proposition~\ref{PropErrorExpansion} and equations \eqref
{EqControlA0}--\eqref{EqControlA2}, \eqref
{EqControlAicentered}, \eqref{EqControlAinoncentered},
\[
{{\mathbb{E}}} \bigl[{\bigl|\DecrAvg{F}^{i;n}(x,q) -{F^{i}}(x,q)\bigr|^2}
\bigr] \leq\frac{ K (\Gamma
^{[\varsigma/2]}_M)^2}{(\Gamma_M)^2} + \frac{K}{\Gamma_M} \leq K'n^{-2(\beta_0\wedge\beta_1)},
\]
implying our claim for $F$, $\DecrAvg{F}^n$. Since $H$ satisfies the
same properties as $F$, the claim follows for $H,\DecrAvg{H}^{n}$.
As a final step, we prove the error control for $\DecrAvg{G}^n$.
Let $\Delta H^n(x,q):=H(x)-\DecrAvg{H}^n(x,q)$ and $E = \{|\Delta
H^n(x,q)|\geq|2H^{-1}|^{-1}\}.$
Markov inequality gives us the control
\[
{ {{\mathbb{P}}}( {E} )} \leq4\bigl |H^{-1}(x)\bigr|^2 {{
\mathbb{E}}} \bigl[{\bigl|\Delta H^n(x,q)\bigr|^2}\bigr] \leq K
n^{-2(\beta_0 \wedge\beta_1)},
\]
which, in conjunction with Theorem~\ref{TheoControlCholesky}, deduces
\begin{eqnarray*}
&&{{\mathbb{E}}} \bigl[{\bigl|G(x)-\DecrAvg{G}^n(x,q)\bigr|^2}
\bigr]
\\
&&\qquad = {{\mathbb{E}}} \bigl[{\bigl|G(x)-\DecrAvg{G}^n(x,q)\bigr|^2
\mathbf{1}_{E}}\bigr] + {{\mathbb{E}}} \bigl[{\bigl|G(x)-
\DecrAvg{G}^n(x,q)\bigr|^2 \mathbf {1}_{E^\complement}}\bigr]
\\
&&\qquad \leq K^\prime n^{-2(\beta_0 \wedge\beta_1)} + {{\mathbb{E}}} \bigl[{\bigl|G(x)-
\DecrAvg {G}^n(x,q)\bigr|^2\mathbf{1}_{E^\complement}}\bigr]
\\
&&\qquad \leq K^{\prime} n^{-2(\beta_0 \wedge\beta_1)} + K n^{-2(\beta_0
\wedge\beta_1)} =
K^{\prime\prime} n^{-2(\beta_0 \wedge\beta_1)}.
\end{eqnarray*}
\upqed\end{pf*}

We can deduce from Proposition~\ref{PropControlDecreasingStep} and the
assumed structure, the following a priori estimates.

\begin{MyCorollary}
\label{CorControlPowerAndDifference}
Under the hypothesis and notation of Proposition~\ref
{PropControlDecreasingStep}, for any $0\leq s\leq T$,
%
\begin{equation}
\label{EqBoundExpectedPower} {{\mathbb{E}}} \bigl[{\bigl|\DecrAvg{F}^n \bigl(
\DecrAvg{X}^n_{\tbar[s]}, \tbar[s] \bigr) \bigr| ^\alpha}
\bigr] \leq K
\end{equation}
and
%
\begin{equation}
\label{EqBoundMeanSquareError} {{\mathbb{E}}} \bigl[{\bigl|\DecrAvg{F}^n \bigl(
\DecrAvg{X}_s^n, \tbar[s] \bigr) - F \bigl(
\DecrAvg{X}_s^n \bigr)\bigr|^2}\bigr] \leq K
n^{-2\beta}.
\end{equation}
The same bounds hold with $\DecrAvg{F}^n, F$ replaced by $\DecrAvg
{H}^n,H$ and $\DecrAvg{G}^n,G$.
\end{MyCorollary}

\begin{pf}
Define
%
\begin{equation}
\label{EqDefinitionFttminusFilter} \bar{\F}_{t,t^-}:= \biggl(\F_t \vee\bigvee
_{q\in\QQ, q < t} \tilde{\F}^q_\infty
\biggr)
\end{equation}
by construction, $\ApproxVar{X}_{\tbar[s]}$ is $\bar{\F}_{\tbar
[s],\tbar
[s]^-}$ measurable and since $\DecrAvg{F}^n(x,\tbar[s])\,\bot\!\!\!
\bot\,\bar{\F
}_{\tbar[s],\tbar[s]^-}$ for any deterministic $x$, we get from
Proposition~\ref{PropControlDecreasingStep},
\[
{{\mathbb{E}}} \bigl[{\bigl|\DecrAvg{F}^n \bigl(\ApproxVar{X}_{\tbar
[s]},
\tbar[s] \bigr)\bigr| ^\alpha}\bigr] ={{\mathbb{E}}} \bigl[{{{\mathbb{E}}}
\bigl[ {\bigl|\DecrAvg{F}^n \bigl( \ApproxVar{X}_{\tbar[s]},\tbar[s]
\bigr)\bigr| ^\alpha} | {\bar{ \F}_{\tbar[s],\tbar[s]^-}} \bigr]}\bigr] \leq{{
\mathbb {E}}} [{K}] = K.
\]
A similar argument leads to \eqref{EqBoundMeanSquareError}, and to the
claims for $\DecrAvg{H}^n,H$ and $\DecrAvg{G}^n,G$.
\end{pf}

Corollary~\ref{CorControlPowerAndDifference} should be understood as
an a priori control on the approximated process. From this control, we
can deduce, using Proposition~\ref{PropAPriori} as in the case of the
effective equation, the existence and strong uniqueness of the solution
of the approximated diffusion \eqref{Eqcontinuousapproximated}. In
addition, Proposition~\ref{PropEstimatePowerApprox} states that
approximation~\eqref{Eqcontinuousapproximated} has a continuous
modification. The result follows from Proposition~\ref
{PropFiniteMomentsTime}, the estimates in Corollary~\ref
{CorControlPowerAndDifference} and Kolmogorov's criterion.

\begin{MyProposition}
\label{PropEstimatePowerApprox}
Under the hypothesis and notation of Proposition~\ref
{PropControlDecreasingStep}, for every $\alpha\geq2$,
\[
{{\mathbb{E}}} \bigl[{\bigl|\ApproxVar{X}^{n}_t -
\ApproxVar{X}^{n}_s\bigr|^\alpha}\bigr] \leq
K_{\alpha
,T} (t-s)^{\alpha/2} \bigl((t-s)^{\alpha/2} +1 \bigr).
\]
Moreover, the solution of \eqref{Eqcontinuousapproximated} has a
continuous modification.
\end{MyProposition}

\subsection{Strong convergence}
\label{SubSecAlmostsurepathwise}

In what follows, we choose $\ApproxVar{X}$ to be continuous in time. We
can proceed to show the mean square convergence of $\ApproxVar{X}^{n}$
toward~$X$.

\begin{MyTheorem}
\label{TheoGeneralConvergence}
Under \ref{HypSlowScale}, \ref{HypFastScale} and
\ref{HypSteps}, let $X$ be defined by \eqref
{EqEffectiveEquation} and $\ApproxVar{X}^{n}$ by \eqref
{Eqcontinuousapproximated}. Let $\beta$ be defined as in \eqref
{EqDefinitionBeta}. Then:
\begin{itemize}
\item if $g\equiv0$ (ODE with random coefficients), then ${{\mathbb
{E}}} [{\sup_{0\leq t \leq T} |X_t-\ApproxVar{X}^{n}_t|^2}] \leq K
n^{-2(1\wedge
\beta)}$;
\item under the full SDE case, ${{\mathbb{E}}} [{\sup_{0\leq t \leq
T} |X_t-\ApproxVar{X}^{n}_t|^2}] \leq K n^{-(1\wedge2\beta)}$.
\end{itemize}
\end{MyTheorem}

\begin{pf} We treat the full SDE case. By definition,
\[
X_{t} - \ApproxVar{X}^{n}_{t} = \int
_0^{t} \bigl[F (X_{s}) - \DecrAvg
{F}^{n} \bigl(\ApproxVar{X}^{n}_{\tbar[s]},\tbar[s]
\bigr) \bigr]\,ds + \int_0^{t} \bigl[ G
(X_{s}) -\DecrAvg{G}^{n} \bigl(\ApproxVar{X}^{n}_{\tbar[s]},
\tbar[s] \bigr) \bigr]\,dW_s.
\]
Our plan is to use Proposition \ref{PropAPriori}(ii). By convexity,
\begin{eqnarray*}
&& \bigl|F (X_{s}) - \DecrAvg{F}^{n} \bigl(
\ApproxVar{X}^{n}_{\tbar
[s]}, \tbar[s] \bigr)\bigr |^2
\\
& &\qquad\leq3 \bigl|F (X_{s}) - F \bigl(\ApproxVar{X}^{n}_s
\bigr) \bigr|^2 + 3\bigl|F \bigl(\ApproxVar{X}^{n}_s
\bigr) - F \bigl(\ApproxVar{X}^{n}_{\tbar[s]} \bigr)
\bigr|^2 + 3\bigl|F \bigl(\ApproxVar{X}^{n}_{\tbar[s]},\tbar[s]
\bigr) - \DecrAvg{F}^{n} \bigl(\ApproxVar{X}^{n}_{\tbar[s]},
\tbar[s] \bigr) \bigr|^2.
\end{eqnarray*}
By Lipschitz assumption in \ref{HypSlowScale},
%
\begin{eqnarray}
{{\mathbb{E}}} \bigl[{\bigl|F (X_{s}) - F \bigl(\ApproxVar{X}^{n}_s
\bigr) \bigr|^2}\bigr] & \leq& K {{\mathbb{E}}} \bigl[{\bigl|X_{s} -
\ApproxVar{X}^{n}_s \bigr|^2}\bigr],
\nonumber
\\[-8pt]
\\[-8pt]
\nonumber
{{\mathbb{E}}} \bigl[{\bigl|F \bigl(\ApproxVar{X}^{n}_s
\bigr) - F \bigl( \ApproxVar{X}^{n}_{\tbar[s]} \bigr)
\bigr|^2}\bigr] & \leq& K {{\mathbb{E}}} \bigl[{\bigl|\ApproxVar{X}^{n}_s
- \ApproxVar{X}^{n}_{\tbar[s]}\bigr|^2}\bigr] \leq K
n^{-1}, \label
{EqEulerError}
\end{eqnarray}
the last inequality being possible for $n$ large enough thanks to
Proposition~\ref{PropEstimatePowerApprox}. Also, by Corollary~\ref
{CorControlPowerAndDifference}, we get
\[
{{\mathbb{E}}} \bigl[{\bigl| F \bigl(\ApproxVar{X}^{n}_{\tbar[s]},
\tbar [s] \bigr) - \DecrAvg{F}^{n} \bigl(\ApproxVar{X}^{n}_{\tbar[s]},
\tbar[s] \bigr) \bigr|^2 }\bigr] \leq K n^{-2\beta}.
\]
Therefore,
%
\begin{equation}
\label{EqBoundF} {{\mathbb{E}}} \bigl[{\bigl|F (X_{s}) -
\DecrAvg{F}^{n} \bigl( \ApproxVar{X}^{n}_{\tbar [s]},
\tbar[s] \bigr) \bigr|^2}\bigr] \leq K \bigl(n^{-(1\wedge2\beta)} + {{
\mathbb{E}}} \bigl[{\bigl|X_s-\ApproxVar {X}^{n}_s\bigr|^2}
\bigr] \bigr).
\end{equation}
Since we may obtain similar bounds for the terms with $G$, we also have
%
\begin{equation}
\label{EqBoundG} {{\mathbb{E}}} \bigl[{\bigl|G (X_{s}) -
\DecrAvg{G}^{n} \bigl( \ApproxVar{X}^{n}_{\tbar [s]},
\tbar[s] \bigr) \bigr|^2}\bigr] \leq K \bigl(n^{-(1\wedge2\beta)} + {{
\mathbb{E}}} \bigl[{\bigl|X_s-\ApproxVar {X}^{n}_s\bigr|^2}
\bigr] \bigr).
\end{equation}
Now, Proposition \ref{PropAPriori}(ii) shows
\[
{{\mathbb{E}}} \bigl[{\bigl|X_{t}-\ApproxVar{X}^{n}_{t}\bigr|^2}
\bigr]  \leq K \int_0^{T}
\bigl(n^{-(1\wedge2\beta)} + {{\mathbb{E}}} \bigl[{\bigl|X_s-
\ApproxVar{X}^{n}_s\bigr|^2}\bigr] \bigr)\,ds.
\]
Therefore, by Gronwall's lemma,
\[
\sup_{0\leq t \leq T} {{\mathbb{E}}} \bigl[{\bigl|X_t-
\ApproxVar{X}^{n}_t\bigr|^2}\bigr] \leq K
n^{-(1\wedge2\beta)}.
\]
Replacing \eqref{EqBoundF} and \eqref{EqBoundG} we get
\[
\sup_{0\leq t \leq T} \bigl( {{\mathbb{E}}} \bigl[{\bigl|F
(X_{s}) - \DecrAvg {F}^{n} \bigl(\ApproxVar{X}^{n}_{\tbar[s]},
\tbar[s] \bigr) \bigr|^2}\bigr] + {{\mathbb{E}}} \bigl[{\bigl|G
(X_{s}) - \DecrAvg{G}^{n} \bigl(\ApproxVar
{X}^{n}_{\tbar [s]},\tbar[s] \bigr) \bigr|^2}\bigr] \bigr)
\\
\leq K n^{-(1\wedge2\beta)}.
\]
So that by Proposition~\ref{PropAPriori},
\[
{{\mathbb{E}}} \Bigl[{\sup_{0\leq t \leq T} \bigl|X_t-
\ApproxVar {X}^{n}_t\bigr|^2}\Bigr] \leq K
n^{-(1\wedge2\beta)}.
\]

Note that the case $g\equiv0$ is proven in the same way, but the Euler
error \eqref{EqEulerError} is bounded by $n^{-2}$ and $G\equiv0$.
This implies the stated result.
\end{pf}

\subsection{Limit distribution}
\label{SubSecLimitDistribution}
In this section we show under slightly stron\-ger regularity assumptions
on the coefficients of the diffusion, that we have convergence in the
weak (uniform topology) sense toward a limit distribution given as the
solution of a particular SDE.

Our plan to prove the limit distribution result is to look at the
rescaled error and its associated stochastic differential equation. We
prove the joint weak convergence of the terms appearing in that SDE and
use the fact that under certain hypothesis the joint convergence of the
terms suffices to deduce the weak convergence of the solution of the
equation. The reader may find most of the needed material on weak
convergence of stochastic integrals and stochastic SDEs in \citet
{jakubowskiconvergence1989}, \citeauthor{kurtzweak1991}
(\citeyear{kurtzweak1991,kurtzweak1996}).

\begin{MyDefinition}
Let $X^{n}$ be a sequence of $\RR^d$-valued semimartingales, and let
$A^{n}(\delta)$ be the predictable process with finite variation null
at zero and $M^{n}(\delta)$ the local martingale null at zero appearing
in the representation of $X^{n}$ as
\[
X^{n}_t = X^{n}_{0} +
A^{n}_t(\delta) + M^{n}_t(\delta) +
\sum_{s\leq
t}\Delta X_s^{n}
\mathbf{1}_{\{|\Delta X_s^n|>\delta\}}.
\]
We say that the sequence $X^{n}$ satisfies property \eqref
{PropertyStar} if for some $\delta> 0$,
{\renewcommand{\theequation}{$\ast$}
\begin{equation}
\bigl\langle M^{n}(\delta),M^{n}(\delta) \bigr
\rangle_T + \int_0^T
\bigl|dA^{n}(\delta)_s\bigr|+ \sum_{s\leq T}
\bigl|\Delta X_s^{n}\bigr|\mathbf{1}_{\{|\Delta X_s^{n}|> \delta\}}
\label{PropertyStar}
\end{equation}}
\hspace*{-3pt}is tight. (The notation $\int_0^T |dA|$ denotes the total variation of
$A$ on $[0,T]$.)
\end{MyDefinition}

The importance of property \eqref{PropertyStar} is shown by the
following theorem; see \citet{jakubowskiconvergence1989}, \citet
{jacodasymptotic1998} and \citet{kurtzweak1996}.

\begin{MyTheorem}
\label{TheoStarConvergenceIntegral}
Let $X^{n}$ be a sequence of $\RR^d$-valued semimartingales relative
to the filtration $\mathcal{F}_t$. Suppose that $X^{n}$ weakly
converges in the Skorokhod topology $D_{\RR^{d_x}}$. Then \eqref
{PropertyStar} is necessary and sufficient for goodness: for any
sequence $H^n$ of $(\F_t)$-adapted c\`{a}dl\`{a}g processes such that
$(H^n,X^n)\Rightarrow(H,X)$ in the Skorokhod topology $D_{M^{d_x\times
d_x}\times\RR^{d_x}}$, then $X$ is a semimartingale w.r.t. the
filtration generated by $(H,X)$ and $(H^n,X^n,\int H^n \,dX^n)
\Rightarrow(H,X,\int H \,dX)$ in the Skorokhod topology $D_{M^{d_x\times
d_x}\times\RR^{d_x}\times\RR^{d_x}}$.
\end{MyTheorem}

Goodness gives us a direct way to show the convergence of sequences of
stochastic integrals, and will play a key role for the convergence of
sequences of SDEs.

Before proceeding to the main propositions of this section, we cite
another useful result concerning weak convergence of sequences of
solutions of SDEs, allowing us to compare the limit of two sequences
with converging coefficients.

\begin{MyTheorem}[{[Theorem~2.5(b) \citet{jacodasymptotic1998}]}]
\label{ThmConvergenceChangeCoefficients}
Consider a sequence of linear SDEs
%
\setcounter{equation}{40}
\begin{equation}
\label{EqLinearSDE} \vartheta_t^n = P_t^n
+\int_0^t \vartheta^n_{s-}
Q_t^n \,dJ_t,
\end{equation}
where the $P_t^n $ are stochastic processes in $\RR^d$, $Q_t^n $ are
stochastic processes in $\RR^{d\times d'}$ and $J_t$ is a
semimartingale in $\RR^{d'}$, and all processes are in same the
filtered probability space. Suppose that we have another sequence of
equations like \eqref{EqLinearSDE} with solution $\vartheta^{\prime
n}$ and coefficients $P^{\prime n}$ and $Q^{\prime n}$. If the sequences
$ \sup_{0\leq s \leq T } \|P^{n}_s\|$ and $\sup_{0\leq s \leq T } \|
Q^{n}_s\| $ are tight, and if
\[
\sup_{0\leq s \leq T }\bigl \|P^{n}_s -
P^{\prime n}_s\bigr\| \mathop{\stackrel{P} {\longrightarrow}}0,\qquad
\sup_{0\leq s \leq T } \bigl\|Q^{n}_s -
Q^{\prime n}_s\bigr\| \mathop{\stackrel{P} {\longrightarrow}}0,
\]
then
\[
\sup_{0\leq s \leq T } \bigl\|\vartheta^{n}_s -
\vartheta^{\prime n}_s\bigr\| \mathop{\stackrel{P} {
\longrightarrow}}0.
\]
\end{MyTheorem}

Proposition~\ref{PropConvergenceOfTuple} shows the weak convergence of
some tuples appearing in the rescaled error SDE.

\begin{MyProposition}
\label{PropConvergenceOfTuple}
Let $\mathcal{I}$ be a set of indices, and consider a family of
independent standard Gaussian variables $ \{\nu^{i;n}_{t_k}\}_{n\in
\NN
^*; 0\leq k \leq n; i\in\mathcal{I}}$ where for any $n,i$ we have
$\nu
^{i;n}_{t_k}$ is $\bar{\F}_{t_k}$ measurable.

Consider the sequence of random processes $A^{0;n}$ (dimension 1),
$A^{1;n},B^{0;n}$ (dimension $d_x$), $B^{2;n}$ (dimension $d_x\times
d_x$), $B^{1;n}$ (dimension $|\mathcal{I}|$) and $B^{3;n}$ (dimension
$|\mathcal{I}|\times d_x$) defined component-wise by
%
\begin{eqnarray}
\label{EqDefHc} B^{0;j;n}_t &:=& \int_0^t
(s-\tbar[s] )\,dW^j_s;\qquad A^{0;n}_t:=
2\int_0^t (s-\tbar[s] )\,ds;
\\
\label{EqDefJc} B_t^{2;l,j;n} &:=& \int_0^t
\sqrt{2} \bigl(W^l_s - W^l_{\tbar[s]}
\bigr) \,d W^j_s;\qquad A^{1;j;n}_t:= \int
_0^t \bigl(W^j_s -
W^j_{\tbar
[s]} \bigr) \,ds;
\\
\label{EqDefKc} B_t^{3;i,j;n} &:=& \int_0^t
\nu^{i;n}_{\tbar[s]} \,dW^j_s;\qquad
B_t^{1;i;n}:= \int_0^t
\nu^{i;n}_{\tbar[s]} \,ds.
\end{eqnarray}
Then we have the following limit results:
%
\begin{eqnarray}
&&\bigl(X,\ApproxVar{X}^{n},W,n A^{0;n}, \sqrt{n}
B^{1;n} \bigr)\quad \Rightarrow \quad\bigl(X,X,W,A^{0},B^{1}
\bigr) \label{EqTupleConv1}
\\
\label{EqTupleConv3}&&
\bigl(X,\ApproxVar{X}^{n}, W,n^{{1}/{2}}A^{0;n},n^{{1}/{2}}B^{0;n},
n^{{1}/{2}}A^{1;n}, n^{{1}/{2}}B^{2;n}_s,B^{1;n},B^{3;n}
\bigr)
\nonumber
\\[-8pt]
\\[-8pt]
\nonumber
&&\qquad\Rightarrow\quad \bigl(X,X,W,0,0,0,B^{2},0,B^{3} \bigr),
\end{eqnarray}
where $A^{0}_t=t$; $B^{0}$, $B^{1}$, $B^{2}$ and $B^{3}$ are standard
Brownian motions defined on an extension of the space $W$, with
dimensions $d_x$, $d_x^2$, $|\mathcal{I}|\times d_x$ and $|\mathcal
{I}|$, respectively.

Moreover, we have $\{ B^{0}$, $B^{2}$, $B^{3}$, $W\}$ are independent;
$\{ B^{0}$, $B^{2}$, $B^{1}$, $W\}$ are independent, and $B^{1;n}$,
$\sqrt{n}B^{2;n}$ and $B^{3;n}$ are ``good'' in the sense of Theorem~\ref{TheoStarConvergenceIntegral}.
\end{MyProposition}

The proof of Proposition~\ref{PropConvergenceOfTuple} will be given in
Section~\ref{SubSecWeakConvTuples}.

\begin{MyProposition}
\label{PropMainDecreasingStepPqCLT} Under the assumptions and notation
of Proposition~\ref{PropControlDecreasingStep}, assume that
$r^y>\varsigma+3$ in \ref{HypSlowScale}, and that there is
$\beta_2 \geq0 $ such
that the asymptotic expansion
%
\begin{equation}
\label{EqConditionGammaMainDecreasingStepLimitDist} \frac{\Gamma
^{[\varsigma/2+1]}_M}{\Gamma_M^{[\varsigma/2]}} = C_2 n^{-\beta_2} + o
\bigl(n^{-\beta_2} \bigr),
\end{equation}
where $\varsigma$ is defined in \eqref{EqDefinitionLbar}, holds. Let
%
\begin{equation}
\label{EqDefinitionRho} \rho= \mathbf{1}_{\{\beta_0>\beta_1\}} \bigl( \beta_2
\wedge(\beta_0-\beta_1) \bigr) + \mathbf{1}_{\{\beta
_0<\beta_1\}}
\bigl(\beta_0 \wedge( \beta_1-\beta_0)
\bigr).
\end{equation}
\begin{longlist}[(ii)]
\item[(i)]
Let $\Phi_F$ be the $d_x\times
d_x$ matrix defined component-wise as
\[
\Phi_F^{i,j}(x) := C_0^{-1}\int
\bigl\langle\sigma^* D_y \phi_{F^i}(x,y), \sigma^*
D_y \phi_{F^j}(x,y) \bigr\rangle\mu^x(dy),
\]
where $\phi_{F^i}$ is the solutions of the Poisson equation \eqref
{EqCenteredPoisson} with source $F^i$. Let
\[
\varphi_F(x):= \mathbf{1}_{\{\beta_1 \geq\beta_0\}} \sqrt{
\Phi_F(x)}; \qquad R_{F}^i(x):=
\mathbf{1}_{\{\beta_0\geq\beta_1\}
} {C_1}\int\bar {v}_{F^i}^{\varsigma,r^y}(x,y)
\mu^x(dy),
\]
with the square root meaning the Cholesky root. Then there exists a
family of independent standard Gaussian variables $ \{\nu^{i;n}_{k}\}
_{n\in\NN^*; 0\leq k \leq n; 1\leq i \leq d_x}$, such that each $\nu
^{i;n}_{k}$ is $\bar{\F}_{t_k}$ measurable and
\[
{{\mathbb{E}}} \Biggl[{\Biggl| n^\beta \bigl( F^i(x)-
\DecrAvg{F}^{i;n}(x,t_k) \bigr) - \sum
_{j=1}^{d_x} \varphi_F^{i,j}(x)
\nu^{j;n}_{k} - R_{F}^i(x)
\Biggr|^2}\Biggr] = O \bigl(n^{-2\rho} \bigr),
\]
for all $x\in\RR^{d_x}$.
\item[(ii)] Under the full SDE case, define in a similar way a $d_x^2$
dimensional random function $R_H$ and a $d_x^2\times d_x^2$ dimensional
random function $\Phi_H$, with
\begin{eqnarray*}
\Phi_H^{i,j,i',j'}(x) &:= &C_0^{-1}\int
\bigl\langle\sigma^* D_y \phi_{H^{i,j}}(x,y), \sigma^*
D_y \phi_{H^{i',j'}}(x,y) \bigr\rangle\mu^x(dy);
\\
\varphi_H(x)&:=& \mathbf{1}_{\{\beta_1\geq\beta_0\}} \sqrt{
\Phi_H(x)}; \\
 R_{H}^{i,j}(x)&:=&
\mathbf{1}_{\{\beta_0 \geq
\beta_1\}} {C_1}\int\bar {v}_{H^{i,j}}^{\varsigma,r^y}(x,y)
\mu^x(dy).
\end{eqnarray*}
Then there exists a family of independent standard Gaussian variables\break 
$\{\nu^{i,j;n}_{k}\}_{n\in\NN^*; 0\leq k \leq n; 0\leq i,j \leq d_x}$,
such that each $\nu^{i,j;n}_{k}$ is $\bar{\F}_{t_k}$ measurable and
\begin{eqnarray*}
&&{{\mathbb{E}}} \Biggl[{\Biggl| n^\beta \bigl( H^{i,j}(x)-\DecrAvg
{H}^{i,j;n}(x,t_k) \bigr) - \sum
_{i',j'=1}^{d_x} \varphi _H^{i,j,i',j'}(x)
\nu^{i',j';n}_{k} - R_{H}^{i,j}(x)
\Biggr|^2}\Biggr] \\
&&\qquad= O \bigl(n^{-2\rho} \bigr),
\end{eqnarray*}
for all $x\in\RR^{d_x}$. Moreover, letting $R_G$, $\varphi_G$ be
defined component-wise for $0\leq i',j'\leq d_x $ as
\begin{eqnarray*}
R_{G}^{i,i} &=& \frac{ R_{H}^{i,i} - 2 \sum_{k=1}^{i-1} R_{G}^{i,k}
G^{i,k}}{2G^{i,i}},\\
 \varphi_{G}^{i,i,i',j'}
&=& \frac{ \varphi_{H}^{i,i,i',j'} - 2 \sum_{j=1}^{i-1} \varphi
_{G}^{i,j,i',j'}G^{i,j}}{2 G^{i,i}},
\end{eqnarray*}
and for $i>j$,
\begin{eqnarray*}
R_{G}^{i,j} &=& \frac{R_{H}^{i,j} -R_{G}^{j,j} G^{i,j} -\sum_{l=1}^{j-1}
(R_{G}^{j,l} G^{i,l}+R_{G}^{i,l}G^{j,l})}{G^{j,j}},
\\
\varphi_{G}^{i,j,i',j'} &=& \frac{ \varphi_{H}^{i,j,i',j'} - \varphi
_{G}^{j,j,i',j'} G^{i,j}-\sum_{l=1}^{j-1} [\varphi
_{G}^{j,l,i',j'} G^{i,l}+\varphi_{G}^{i,l,i',j'} G^{j,l}] }{G^{j,j}}.
\end{eqnarray*}
\end{longlist}
Then
\begin{eqnarray*}
&&{{\mathbb{E}}} \Biggl[{\Biggl| n^\beta \bigl( G^{i,j}(x)-\DecrAvg
{G}^{i,j;n}(x,t_k) \bigr) - \sum
_{i',j' =1}^{d_x} \varphi _G^{i,j,i',j'}(x)
\nu^{i',j';n}_{k} - R_G^{i,j}(x)
\Biggr|^2}\Biggr]\\
&&\qquad = O \bigl(n^{-2\rho} \bigr).
\end{eqnarray*}
\end{MyProposition}

\begin{pf}
(i) We prove the first claim. We use the expansion of Proposition~\ref
{PropErrorExpansion} up to order $\varsigma$ as in Proposition~\ref
{PropControlDecreasingStep}, and estimates
\eqref{EqControlA0}--\eqref{EqControlA2},
\eqref{EqControlAicentered}--\eqref{EqControlZi} to get for any
$x$ that
%
\begin{eqnarray}
\label{EqMainErrorTerms} %
&&{{\mathbb{E}}} \biggl[{\biggl| \bigl({F^i}(x)-
\DecrAvg{{F}}^{i;n}(x,q) \bigr) - \frac{1}{{\Gamma_M}}
\bigl(N_{{F^i},M}(x,q)+A_{{F^i},M}^{(\varsigma
)}(x,q) \bigr)
\biggr|^2}\biggr]
\nonumber
\\[-8pt]
\\[-8pt]
\nonumber
&&\qquad= O \bigl( (\Gamma_M)^{-2} \bigl[1 + \bigl(
\Gamma_M^{[\varsigma
/2+1]} \bigr)^2 \bigr] \bigr).
\end{eqnarray}

Let us examine separately three cases depending on the relation between
$\beta_0$ and~$\beta_1$:
\begin{itemize}
\item If $\beta_0>\beta_1$: In this case $\beta=\beta_1$, and by
definition of $\beta_1$ it follows that
%
\begin{eqnarray}\label{EqExpansionLimitBigBeta0}
&&{{\mathbb{E}}} \bigl[{\bigl| n^\beta \bigl( F^i(x)-
\DecrAvg{F}^{i;n}(x,q) \bigr) - R_{F}^i(x)
\bigr|^2}\bigr]
\nonumber\\
&&\qquad \leq K {{\mathbb{E}}} \bigl[{\bigl| \bigl(\Gamma^{[\varsigma/2]}
\bigr)^{-1} \Gamma_M \bigl( F^i(x)-
\DecrAvg{F}^{i;n}(x,q) \bigr) - R_{F}^i(x)
\bigr|^2}\bigr]
\nonumber
\\
&&\qquad \leq K' {{\mathbb{E}}} \bigl[\bigl|\Gamma^{[\varsigma/2]}_M\bigr|^{-2}
\bigl|\Gamma _M \bigl(F^i(x)-\DecrAvg{F}^{i;n}(x,q)
\bigr)
\nonumber
\\[-8pt]
\\[-8pt]
\nonumber
&&\hspace*{105pt}{}- N_{F^i,M}(x,q)- A_{F^i,M}^{\varsigma }(x,q)
\bigr|^2\bigr]
\\
&&\qquad\quad{}
+ K' {{\mathbb{E}}} \bigl[{\bigl| \bigl( \Gamma^{|[\varsigma/2]}_M
\bigr)^{-1} \bigl(A_{F^i,M}^{\varsigma}(x,q) -
R_{F}^i(x) \bigr)\bigr|^2 }\bigr]\nonumber \\
&&\qquad\quad{}+
K' {{\mathbb{E}}} \bigl[{\bigl| \bigl( \Gamma^{|[\varsigma/2]}_M
\bigr)^{-1} N_{F^i,M}(x,q)\bigr |^2 }\bigr].
\nonumber
\end{eqnarray}

The first term in the right-hand side of \eqref
{EqExpansionLimitBigBeta0} can be controlled by rescaling~\eqref
{EqMainErrorTerms} to get
%
\begin{eqnarray}\label{EqErrorPsiB0greaterB1}
&&{{\mathbb{E}}} \bigl[{\bigl|\Gamma^{[\varsigma/2]}_M\bigr|^{-2} \bigl|
\Gamma_M \bigl(F^i(x)-\DecrAvg{F}^{i;n}(x,q)
\bigr) - N_{F^i,M}(x,q)- A_{F^i,M}^{\varsigma}(x,q)
\bigr|^2}\bigr]
\nonumber
\\[-8pt]
\\[-8pt]
\nonumber
&&\qquad= O \bigl( \bigl(\Gamma_M^{[\varsigma/2]} \bigr)^{-2}
\bigl[ \bigl(\Gamma_M^{[\varsigma/2+1]} \bigr)^2 + 1 \bigr]
\bigr).
\end{eqnarray}
From \eqref{EqControlN} we control the third term in the right-hand
side of \eqref{EqExpansionLimitBigBeta0}
%
\begin{equation}
\label{EqOrderNPsi} {{\mathbb{E}}} \bigl[{\bigl| \bigl(\Gamma_M^{[\varsigma
/2]}
\bigr)^{-1} N_{F^i,M}(x,q)\bigr|^2}\bigr] = O \bigl(
\bigl( \Gamma_M^{[\varsigma/2]} \bigr)^{-2}
\Gamma_M \bigr).
\end{equation}
To control the second term of \eqref{EqExpansionLimitBigBeta0}, let us
define
%
\begin{equation}
\label{EqDefinitionAbar} \bar{A}_{F^i,M}^{\varsigma}(x,q):= \sum
_{k=1}^{M} \gamma_k^{\varsigma
/2}
\bar{v}^{\varsigma,r^y}_{F^i} \bigl(x,\bar{Y}^{x,q}_{k-1}
\bigr)
\end{equation}
for $\bar{v}^{\varsigma,r^y}_{F^i}$ defined in \eqref
{EqDefinitionvbar}. We can compare $(\Gamma^{[\varsigma/2]}_M)^{-1}
A^{\varsigma}_{F^i,M}$ and $(\Gamma^{[\varsigma/2]}_M)^{-1} \times  \bar
{A}^{\varsigma}_{F^i,M}$ in $L_2$ by \eqref{EqEstimateSquareCentered}.
Indeed, thanks to controls \eqref{EqControlAicentered} and
\eqref{EqControlAinoncentered}, and the fact that for some $K\in\RR
^+$, $\Gamma_M^{[\varsigma]} \leq K \Gamma_M$, we have
%
\begin{eqnarray}\label{EqAminusAbar}
&& {{\mathbb{E}}} \bigl[{\bigl| \bigl(\Gamma^{[\varsigma/2]}_M
\bigr)^{-1} \bigl( A^{\varsigma }_{F^i,M}(x,q) -
\bar{A}^{\varsigma
}_{F^i,M}(x,q) \bigr) \bigr|^2}\bigr]\nonumber
\\
&&\qquad = {{\mathbb{E}}} \Biggl[{\Biggl| \bigl(\Gamma^{[\varsigma/2]}_M
\bigr)^{-1}\sum_{k=1}^{M}
\gamma_k^{\varsigma/2} \bigl( v_{F^i}^{\varsigma,r^y}
\bigl(x,\bar{Y}^{x,q}_{k-1},U_k^q \bigr)
- \bar{v}^{\varsigma,r^y}_{F^i} \bigl(x,\bar{Y}^{x,q}_{k-1}
\bigr) \bigr)\Biggr |^2}\Biggr]
\\
&&\qquad = O \bigl( \bigl(\Gamma^{[\varsigma/2]} \bigr)^{-2}
\Gamma_M \bigr).
\nonumber
\end{eqnarray}

It remains to show that
%
\begin{equation}
\label{EqAbarminusR}{{\mathbb{E}}} \bigl[{\bigl| \bigl(\Gamma^{[\varsigma
/2]}_M
\bigr)^{-1}\bar{A}_{F^i,M}^{\varsigma}(x,q) +
R_F^i(x) \bigr|^2}\bigr]=O \bigl(n^{-2\rho}
\bigr).
\end{equation}
Indeed, from the definition of $\beta_0$ and $\beta_1$, $\Gamma
^{[\varsigma/2]}_M = O(n^{\beta_0-\beta_1})$ so that it diverges.
Moreover, from the assumed regularity hypothesis, $\bar{v}^{\varsigma
,r^y}_{F^i}(x,y)$ is $C_{p,b}^{2,r^y-\varsigma}$. Therefore,
Proposition~\ref{PropPropertiesPhi} guarantees the existence of a
solution to the centered Poisson equation with source $\bar
{A}^{\varsigma,r^y}_{F^i}(x,y)$ of the same regularity, and thus
Proposition~\ref{PropControlDecreasingDifferentStep} shows that $\bar
{A}^{\varsigma}_{F^i,M}(x,q)$ converges uniformly with respect to $x$
in $L^2$ to $-R_F^i(x)$ with rate $(\beta_0-\beta_1) \wedge\beta_2
\geq\rho$ since
\begin{eqnarray*}
&& \bigl(\Gamma_M^{\varsigma/2} \bigr)^{-2} \bigl(1+
\Gamma _M^{[\varsigma]}+\Gamma_M^{[\varsigma-1]}+ \bigl(
\Gamma_M^{[\varsigma
/2+1]} \bigr)^2 \bigr)
\\
&&\qquad \leq K \bigl(\Gamma_M^{\varsigma/2} \bigr)^{-2}
\bigl( \Gamma_M+ \bigl(\Gamma_M^{[\varsigma/2+1]}
\bigr)^2 \bigr)= O \bigl(n^{-2((\beta_0-\beta_1 )\wedge\beta
_2)} \bigr).
\end{eqnarray*}
The claim follows from replacing \eqref{EqErrorPsiB0greaterB1}, \eqref
{EqOrderNPsi}, \eqref{EqAminusAbar} and \eqref{EqAbarminusR} in
\eqref{EqExpansionLimitBigBeta0}.
\item If $\beta_1>\beta_0$,
we follow a similar approach. We expand the rescaled error term to find
%
\begin{eqnarray}\label{EqExpansionLimitBigBeta1}
&&{{\mathbb{E}}} \bigl[{\bigl| n^\beta \bigl( F^i(x)-
\DecrAvg{F}^{i;n}(x,q) \bigr) - R_{F}^i(x)
\bigr|^2}\bigr] \nonumber
\\
&&\qquad \leq K' {{\mathbb{E}}} \bigl[|\Gamma_M|^{-1/2}
\bigl| \Gamma_M \bigl(F^i(x)-\DecrAvg{F}^{i;n}(x,q)
\bigr) \nonumber\\
&&\hspace*{40pt}\qquad\quad{}- N_{F^i,M}(x,q)- A_{F^i,M}^{\varsigma}(x,q)
\bigr|^2\bigr]
\\
&&\qquad\quad{} + K' {{\mathbb{E}}} \bigl[{\bigl|( \Gamma_M)^{-1/2}
A_{F^i,M}^{\varsigma
}(x,q)\bigr |^2 }\bigr]\nonumber
\\
&&\quad\qquad{} + K' {{\mathbb{E}}} \bigl[{\bigl|( \Gamma_M)^{-1/2}
\bigl(N_{F^i,M}(x,q) - \Phi_F \bigr)\bigr|^2 }\bigr].
\nonumber
\end{eqnarray}

By rescaling \eqref{EqMainErrorTerms} we get
\begin{eqnarray*}
&&{{\mathbb{E}}} \bigl[{|\Gamma_M|^{-1/2}\bigl |
\Gamma_M \bigl(F^i(x)-\DecrAvg {F}^{i;n}(x,q)
\bigr) - N_{F^i,M}(x,q)- A_{F^i,M}^{\varsigma}(x,q)
\bigr|^2}\bigr]
\\
&&\qquad= O \bigl( (\Gamma_M)^{-1} \bigl[ 1+ \bigl(
\Gamma_M^{[\varsigma
/2+1]} \bigr)^2 \bigr] \bigr),
\end{eqnarray*}
and from \eqref{EqControlAinoncentered},
\[
{{\mathbb{E}}} \bigl[{\bigl|(\Gamma_M)^{-1/2}
A_{F^i,M}^{\varsigma}(x,q) \bigr|^2}\bigr] = O \bigl( (
\Gamma_M)^{-1} \bigl( \Gamma_M^{[\varsigma/2]}
\bigr)^2 \bigr).
\]
So it remains to consider the $N_M$ term. Note that since the $U_k^q$
are independent standard Gaussian vectors, $(C_0 \sqrt{\Gamma_M})^{-1}
N_{F^i,M}(x,q)$ when $i$ ranges $1,\ldots,d_x$ is a Gaussian vector.

Let us study its covariance matrix $\Phi_F^{n}$. Using \eqref
{EqEstimateProductCentered} we get for $i,j=1,\ldots,n$
\begin{eqnarray*}
\Phi_F^{i,j;n}(x,q) &:=& {{\mathbb{E}}} \biggl[{
\frac{1}{\Gamma_M} N_{F^i,M}(x,q) N_{F^j,M} \bigl(x,q'
\bigr) }\biggr]
\\
& =&\mathbf{1}_{\{q=q'\}}\sum_{k,k'=1}^{M}
\gamma_k \bigl\langle\sigma^*(\cdot)D_y\phi
_{F^i}(\cdot), \sigma^*(\cdot)D_y\phi_{F^{j}}(
\cdot) \bigr\rangle \bigl(x,\DecrAvg{Y}^{x,q}_{k-1} \bigr).
\end{eqnarray*}

Define $\varphi^n_F=\sqrt{\Phi^n_F}$ (the Cholesky decomposition).
Then, there exists a family of independent Gaussian variables $\nu
_{t_k}^{i,j;n}$, $\bar{\F}_{t_k}$-measurable such that
\[
(\Gamma_M)^{-1} N_{F^i,M}(x,q) = \sum
_{j=1}^{d_x} \varphi^{i,j;n}\nu
_{t_k}^{i,j;n}.
\]
Moreover, from Proposition~\ref{PropConvergenceEstimator} and
Proposition~\ref{PropControlDecreasingStep},
we have that $\Phi^n_F(x,q) $ converges uniformly in $x$ in $L^2$ to
$\Phi_F(x)$ as defined in the claim with rate $O(n^{-\beta})$. By
Theorem~\ref{TheoControlCholesky} we get the same uniform convergence
for $\varphi_F^n$. The claim follows in this case.
\item The case $\beta_0 = \beta_1$ is straightforward from what has
been proven in the previous cases.
\end{itemize}

(ii) Since $H, \DecrAvg{H}^n$ satisfy the same properties as
$F,\DecrAvg{F}^n$, we get the claim for $R_H$, $\varphi_H$ and $\nu
_k^{i,j;n}$ by analogous arguments. Replacing this result in the
sensitivity of the Cholesky procedure given in Lemma~\ref
{LemErrorCholesky}, and taking into account the independence of the
Gaussian entries, we get the claim for $R_G$ and $\varphi_G$.
\end{pf}

Let $\{\upsilon_n\}$ be a sequence of increasing positive numbers, and
let us consider the sequence of rescaled error processes $\zeta^{n}$,
defined by
\[
\zeta_t^{n}:= \upsilon_n
\bigl(X_t - \ApproxVar{X}^{n}_t \bigr).
\]

We can show that this sequence of processes converges in distribution
in the uniform convergence topology to a process $\zeta$ defined as the
solution to a certain stochastic differential equation. We divide the
analysis in two main cases: a first one in which $G(x)\equiv0 $, that
is, when $X$ is the solution to an ordinary differential equation, and
the case when $G(x)$ is nondegenerate. Just as in the asymptotic error
obtained for the usual stochastic Euler method given in \citet
{jacodasymptotic1998}, we will obtain different rates and different
components in the equation for both cases.

\begin{MyTheorem}[(Limit distribution)]
\label{TheoAbstractLimitDistribution}
Under the assumptions and notation of Proposition~\ref
{PropMainDecreasingStepPqCLT}, let $\rho,R_F,\varphi_F,R_G,\varphi_G$
be defined as in Proposition~\ref{PropMainDecreasingStepPqCLT} and
$\beta$ defined in \eqref{EqDefinitionBeta}.
\begin{longlist}[(i)]
\item[(i)] 
[ODE case-$G(x)\equiv
0$.] Let $B^{1}$ be the Brownian process given in Proposition~\ref
{PropConvergenceOfTuple}.
Let $r = 1 \wedge(1/2 + \beta)$, and suppose $\rho\geq r-\beta$. Let
\[
\zeta_t^{n}:= n^r \bigl(X_t -
\ApproxVar{X}^{n}_t \bigr).\
\]
Then $\zeta^n\Rightarrow\zeta^\infty$ in the uniform convergence
sense, where $\zeta^\infty$ is solution of the system
%
\begin{eqnarray}\label{EqLimitDistODE}
\zeta_t^{\infty,i} &=& \sum_{j=1}^{d_x}
\biggl( \int_0^t \partial_{x^j}
F^i(X_s)\zeta^{\infty,j}_s \,ds +
\mathbf{1}_{\{\beta\geq1/2\}} \frac{1}{2} \int_0^t
\partial_{x^j} F^i(X_s) F^j(X_s)
\,ds \biggr)
\nonumber
\\[-8pt]
\\[-8pt]
\nonumber
&&{}+ \mathbf{1}_{\{\beta\leq1/2\}} \Biggl( \int_0^t
R_{F}^{i}(X_s) \,ds + \sum
_{l=1}^{d_x}\int_0^t
\varphi_F^{i,l}(X_s) \,dB^{1;l}_s
\Biggr).
\end{eqnarray}
\item[(ii)] 
[SDE case-$G(x)\neq
0$.] Let $B^{2}$ and $B^{3}$ be the independent Brownian processes given
in Proposition~\ref{PropConvergenceOfTuple}.
Let $r = (1/2 \wedge\beta)$ and
\[
\zeta_t^{n}:= n^r \bigl(X_t -
\ApproxVar{X}^{n}_t \bigr).
\]
Then $\zeta^n \Rightarrow\zeta^\infty$, where $\zeta^\infty$ is
solution of the system for $i=1,\ldots,d_x$ of
%
\begin{eqnarray}
\label{EqlimiDistSDE} %
 \zeta_t^{\infty,i} &=&
\sum_{j} \biggl( \int_0^t
\partial_{x^j} F^i(X_s)\zeta^{\infty,j}_s
\,ds + \int_0^t R_G^{i,j}(X_s)
\,dW^j_s \biggr)\nonumber \\
&&{}+ \mathbf{1}_{\{\beta\leq1/2\}} \sum
_{j,k,l =1}^{d_x}\int_0^t
\varphi_G^{i,j,l,k}(X_s) \,dB^{3;l,k,j}_s
\nonumber
\\[-8pt]
\\[-8pt]
\nonumber
&&{} + \mathbf{1}_{\{\beta\leq1/2\}} \sum_{j,l=1}^{d_x}
\int_0^t \partial_{x^j}
G^{i,l}(X_s)\zeta^{\infty,j}_s
\,dW^l_s \\
&&{}+ \mathbf{1}_{\{\beta\geq1/2\}} \frac
{1}{\sqrt{2}}
\sum_{j,k,l=1}^{d_x} \int_0^t
\partial_{x^j} G^{i,l}(X_s)
G^{j,k}(X_s) \,dB^{2;k,l}_s.\nonumber
\end{eqnarray}
\end{longlist}
\end{MyTheorem}

Let us remark that if $\beta> 1/2$ in Theorem~\ref
{TheoAbstractLimitDistribution}, the error of the Euler scheme
dominates: we recover the limit distribution error for an Euler scheme
with exact coefficients given in \citet{kurtzwong-zakai1991} or \citet
{jacodasymptotic1998}. By contrast, if $\beta<1/2$, it is the
decreasing Euler estimate error that becomes dominant. Since a higher
$\beta$ is generally only achieved by paying a higher price in the
required number of steps for the decreasing Euler step, the optimal
choice implies fixing $\beta=1/2$.

Before proving Theorem~\ref{TheoAbstractLimitDistribution}, let us
show how it implies Theorem~\ref{TheoShortMainDecreasingStep}.

\begin{pf*}{Proof of Theorem~\ref{TheoShortMainDecreasingStep}}
The result is obtained, from Theorems \ref{TheoGeneralConvergence}
and~\ref{TheoAbstractLimitDistribution}, since \ref{HypSlowScale}
and \ref{HypFastScale} are
directly assumed and as the sequence defined as $\gamma_k = \gamma_1
k^{-\theta}$ for $0<\theta<1$ satisfies Hypothesis
\ref{HypSteps}. Moreover,
recall that we fixed $M(n)= \lceil M_1n^{1/(1-\theta)} \rceil$, and we
have for $n$ large enough,
\begin{eqnarray*}
\Gamma_M &\approx&\frac{\gamma_0 M_1^{1-\theta}n}{1-\theta},\qquad \frac
{\Gamma_M^{[\varsigma/2]}}{\Gamma_M} \approx
\frac{(1-\theta)
M_1^{-(\varsigma/2-1)\theta}
 n^{-{(\varsigma/2-1)\theta
}/{(1-\theta)
}}}{1-\varsigma\theta/2},
\\
\frac{\Gamma_M^{[\varsigma/2+1]}}{\Gamma_M^{[\varsigma/2]}} &\approx&\frac{ (1-\varsigma\theta/2)
M_1^{-\theta} n^{-
{\theta
}/{(1-\theta)
}}}{1-(\varsigma/2+1)\theta},
\end{eqnarray*}
so that we get from Proposition~\ref{PropMainDecreasingStepPqCLT},
that $\beta_0 = 1/2$ and
\[
\beta_1=\frac{(\varsigma/2-1)\theta}{1-\theta},\qquad \beta_2=\frac
{\theta
}{1-\theta},\qquad
C_0 \approx\frac{\gamma_0 M_1^{1-\theta}}{1-\theta},\qquad C_1 \approx
\frac{(1-\theta) M_1^{-\theta}}{1-2\theta}.
\]

Recall that $\varsigma$ is defined in \eqref{EqDefinitionLbar} and
stands for the first nonzero term in the error expansion of the
decreasing Euler estimator. Let us assume we are in the worst case when
it attains its minimal value $\varsigma=4$. Hence
\[
\beta_1 = \frac{\theta}{1-\theta},\qquad C_1 =
\frac{(1-\theta)
M_1^{-\theta
}}{1-2\theta}.
\]
Let us now deduce
the conditions on $\theta$ are then deduced from the conditions in
Theorem~\ref{TheoAbstractLimitDistribution} for each of our study cases:
\begin{itemize}
\item ODE with random coefficients: From the conditions of Theorem~\ref
{TheoAbstractLimitDistribution} we have
\[
r = 1 \wedge \bigl(\tfrac{1}{2}+\beta \bigr) = \tfrac{1}{2} +
\bigl( \tfrac
{1}{2}\wedge\beta \bigr) = \tfrac{1}{2} + \bigl(
\tfrac{1}{2}\wedge(\beta_0 \wedge\beta_1) \bigr) =
\tfrac{1}{2} + \beta
\]
since we should verify $\rho\geq r-\beta= 1/2 $, this implies
\[
|\beta_0 -\beta_1| = \biggl|\frac{1}{2}-
\frac{\theta}{1-\theta}\biggr |\geq\rho\geq\frac{1}{2},
\]
which is the case if $\theta\in[1/2,1)$. Moreover, since in this case
$\beta_1\geq1>\beta_0=1/2$, we get $r=1/2$, and the $R_F$ term disappears.
\item Full SDE case: We have $r =\beta= 1/2 \wedge(\theta/(1-\theta
))$ the only restriction comes from imposing $\beta=1/2$. This is
obtained for $1/3\leq\theta<1$. Note that the $R_G$ term is different
from zero only if $\theta=1/3$.
\end{itemize}
Finally, note that if $\varsigma> 4$, we get from the constraints
$\theta\in[1/2,1)$ in the ODE with random coefficients case that
$\beta
_1> \beta_0+1/2$ and from fixing $\theta\in[1/3,1)$ in the full
SDE case that $\beta_1>\beta_0=1/2$, $\beta_1>\beta_2$. In both those
cases the $R_G$ term is zero.
\end{pf*}

\begin{MyRemark}
It should be noted from the proof of Theorem~\ref
{TheoShortMainDecreasingStep} that knowing a priori that $\varsigma>4$
makes it possible to obtain a lower inferior bound for $\theta$ in the
theorem. Since in general we do not know $\varsigma$, we have stated
our results with the sometimes sub-optimal limits.
\end{MyRemark}

\begin{pf*}{Proof of Theorem~\ref{TheoAbstractLimitDistribution}}
(a) Let us deal first with the full SDE case. We have from the
definition of $\zeta^{n}$ that
%
\begin{equation}\qquad
\zeta^{n}_t = \int_0^t
n^{r} \bigl( F(X_s) - \DecrAvg{F}^{n} \bigl(
\ApproxVar{X}^{n}_{\tbar[s]},\tbar[s] \bigr) \bigr) \,ds + \int
_0^t n^{r} \bigl(
G(X_s) - \DecrAvg{G}^{n} \bigl(\ApproxVar{X}^{n}_{\tbar[s]},
\tbar[s] \bigr) \bigr) \,dW_s.\label{EqdevelopmentZita}
\end{equation}
Let us examine each one of these terms separately. Denoting by $x^i$
the $i$th component of $x$, let $x,y \in\RR^{d_x}$. We define the set
of vectors $\Delta^j (x,y)$
\[
\Delta^j(x,y):= \cases{ x, &\quad $\mbox{for } j=0,$\vspace*{2pt}
\cr
(y_1,y_2,\ldots,y_j,x_{j+1},x_{j+2},
\ldots,x_{d_x})^*, &\quad $\mbox{for } 1\leq j\leq d_x$,}
\]
and
\[
\Delta^j F^i (x,y):= \mathbf{1}_{\{x^j \neq y^j \}} \biggl(
\frac{
F^i(\Delta^{j-1}(x,y)) - F^i(\Delta^j(x,y)) }{
x^j - y^j } \biggr) + \mathbf{1}_{\{x^j = y^j \}} \partial_{x^j}F^i(x),
\]
and recalling that
\[
\ApproxVar{X}^{j,n}_{s} - \ApproxVar{X}^{j,n}_{\tbar[s]}
= F^j \bigl(\ApproxVar{X}^{n}_{\tbar[s]} \bigr)
\bigl(s-\tbar[s] \bigr) + \sum_{l=1}^{d_x}
G^{j,l} \bigl(\ApproxVar{X}^{n}_{\tbar[s]} \bigr)
\bigl(W^l_s-W^l_{\tbar[s]} \bigr),
\]
we have
\begin{eqnarray*}
&&\int_0^t n^{r} \bigl[
F^i(X_s) - \DecrAvg{F}^{i;n} \bigl(\ApproxVar
{X}^{n}_{\tbar[s]},\tbar[s] \bigr) \bigr] \,ds \\
&&\qquad= \int
_0^t n^{r} \bigl(
F^i(X_s) - F^i \bigl(\ApproxVar{X}^{n}_{s}
\bigr) \bigr) \,ds + \int_0^t n^{r}
\bigl( F^i \bigl(\ApproxVar{X}^{n}_{s} \bigr) -
F^i \bigl(\ApproxVar{X}^{n}_{\tbar[s]} \bigr) \bigr)
\,ds
\\
&&\qquad\quad{}+ \int_0^t n^{r} \bigl(
F^i \bigl(\ApproxVar{X}^{n}_{\tbar[s]} \bigr) -
\DecrAvg{F}^{i;n} \bigl(\ApproxVar{X}^{n}_{\tbar[s]},
\tbar[s] \bigr) \bigr) \,ds
\end{eqnarray*}
so that
\begin{eqnarray*}
&&\int_0^t n^{r} \bigl[
F^i(X_s) - \DecrAvg{F}^{i;n} \bigl(\ApproxVar
{X}^{n}_{\tbar[s]},\tbar[s] \bigr) \bigr] \,ds \\
&&\qquad = \int
_0^t \sum_j
\Biggl[  n^{r} \Delta^j F^i
\bigl(X_s, \ApproxVar{X}^{n}_{s} \bigr)
\bigl(X^j_s - \ApproxVar{X}^{j;n}_{s}
\bigr)
\\
&&\hspace*{33pt}\qquad\quad{} +n^{r}\Delta^j F^i \bigl(
\ApproxVar{X}^{n}_s, \ApproxVar{X}^{n}_{\tbar
[s]}
\bigr) F^j \bigl(\ApproxVar{X}^{n}_{\tbar[s]} \bigr)
(s-\tbar[s] )
\\
&&\hspace*{33pt}\qquad\quad{} +\sum_{l=1}^{d_x}n^{r}
\Delta^j F^i \bigl(\ApproxVar{X}^{n}_s,
\ApproxVar{X}^{n}_{\tbar[s]} \bigr) G^{j,l} \bigl(
\ApproxVar{X}^{n}_{\tbar[s]} \bigr) \bigl(W^l_s-W^l_{\tbar[s]}
\bigr) \Biggr]\,ds
\\
&&\qquad\quad{} + n^{r-\beta}\int_0^t
n^{\beta} \bigl( F^i \bigl(\ApproxVar{X}^{n}_{\tbar
[s]}
\bigr) - \DecrAvg{F}^{i;n} \bigl(\ApproxVar{X}^{n}_{\tbar[s]},
\tbar[s] \bigr) \bigr) \,ds.
\end{eqnarray*}
Following the same approach we obtain for each $l=1,\ldots,d_x$,
\begin{eqnarray*}
&&\int_0^t  n^{r} \bigl[
G^{i,l}(X_s) - \DecrAvg{G}^{i,l;n} \bigl(
\ApproxVar{X}^{n}_{\tbar[s]},\tbar[s] \bigr) \bigr]
\,dW^l_s
\\
&&\qquad = \int_0^t \sum
_j \Biggl[  n^{r}
\Delta^j G^{i,l} \bigl(X_s,
\ApproxVar{X}^{n}_{s} \bigr) \bigl(X^j_s
- \ApproxVar{X}^{j;n}_{s} \bigr)
\nonumber
\\
&&\hspace*{32pt}\qquad\quad{} +n^{r}\Delta^j G^{i,l} \bigl(
\ApproxVar{X}^{n}_s, \ApproxVar{X}^{n}_{\tbar
[s]}
\bigr) F^j \bigl(\ApproxVar{X}^{n}_{\tbar[s]} \bigr)
(s-\tbar[s] )
\nonumber
\\
&&\hspace*{32pt}\qquad\quad{}+\sum_{k=1}^{d_x}n^{r}
\Delta^j G^{i,l} \bigl(\ApproxVar{X}^{n}_s,
\ApproxVar{X}^{n}_{\tbar[s]} \bigr) G^{j,k} \bigl(
\ApproxVar{X}^{n}_{\tbar[s]} \bigr) \bigl(W^k_s-W^k_{\tbar[s]}
\bigr) \Biggr]\,dW^l_s
\\
&&\qquad\quad{} + n^{r-\beta}\int_0^t
n^{\beta} \bigl( G^{i,l} \bigl(\ApproxVar{X}^{n}_{\tbar[s]}
\bigr) - \DecrAvg{G}^{i,l;n} \bigl(\ApproxVar{X}^{n}_{\tbar[s]},
\tbar[s] \bigr) \bigr) \,dW^l_s.
\nonumber
\end{eqnarray*}
By identifying terms in the obvious way, we write
\[
\zeta_t^{i,n} = \bigl( P^{i,n}_1(t)+P^{i,n}_2(t)
\bigr) + \int_0^t \bigl\langle
Q^{i;n}_1(s), \zeta_s^{n} \bigr
\rangle \,ds + \sum_{l=1}^{d_x}\int
_0^t \bigl\langle Q^{i,l;n}_2(s),
\zeta_s^{n} \bigr\rangle \,dW^l_s,
\]
where $Q_1^i$, $Q_2^{i,l}$ are $d_x$ dimensional random processes with
components
\[
Q^{j,i;n}_1(s) = \Delta^j F^i
\bigl(X_s, \ApproxVar{X}^{n}_{s} \bigr),\qquad
Q^{j,i,l;n}_2(s) = \Delta^j G^{i,l}
\bigl(X_s, \ApproxVar{X}^{n}_{s} \bigr)
\]
and
\begin{eqnarray*}
P^{i;n}_2(s) &= & n^{r-\beta}\int
_0^t n^{\beta} \bigl( F^i
\bigl(\ApproxVar{X}^{n}_{\tbar[s]} \bigr) - \DecrAvg{F}^{i;n}
\bigl(\ApproxVar{X}^{n}_{\tbar
[s]},\tbar[s] \bigr) \bigr)\,ds
\\
& &{}+ n^{r-\beta}\int_0^t
n^{\beta} \bigl( G^{i,l} \bigl(\ApproxVar{X}^{n}_{\tbar
[s]}
\bigr) - \DecrAvg{G}^{i,l;n} \bigl(\ApproxVar{X}^{n}_{\tbar[s]},
\tbar[s] \bigr) \bigr) \,dW^l_s.\vspace*{-9pt}
\end{eqnarray*}
\begin{eqnarray*}
 P^{i;n}_1(s) &=& \int_0^t
\sum_j \Biggl[n^{r}
\Delta^j F^i \bigl(\ApproxVar{X}^{n}_s,
\ApproxVar {X}^{n}_{\tbar
[s]} \bigr) F^j \bigl(
\ApproxVar{X}^{n}_{\tbar[s]} \bigr) (s-\tbar[s] )\\
&&\hspace*{34pt}{} +
\sum_{l=1}^{d_x}n^{r}
\Delta^j F^i \bigl(\ApproxVar{X}^{n}_s,
\ApproxVar{X}^{n}_{\tbar[s]} \bigr) G^{j,l} \bigl(
\ApproxVar{X}^{n}_{\tbar[s]} \bigr) \bigl(W^l_s-W^l_{\tbar[s]}
\bigr) \Biggr]\,ds
\\
&&{} + \int_0^t \sum
_{j,l=1}^{d_x} \Biggl[ n^{r}
\Delta^j G^{i,l} \bigl(\ApproxVar{X}^{n}_s,
\ApproxVar{X}^{n}_{\tbar
[s]} \bigr) F^j \bigl(
\ApproxVar{X}^{n}_{\tbar[s]} \bigr) (s-\tbar[s] ) \\
&&\hspace*{51pt}{}\times\sum
_{k=1}^{d_x}n^{r}
\Delta^j G^{i,l} \bigl(\ApproxVar{X}^{n}_s,
\ApproxVar{X}^{n}_{\tbar[s]} \bigr) G^{j,k} \bigl(
\ApproxVar{X}^{n}_{\tbar[s]} \bigr) \bigl(W^k_s-W^k_{\tbar[s]}
\bigr) \Biggr] \,dW^l_s.
\end{eqnarray*}

(b) In this step, we introduce a nicer diffusion and study its
convergence, and prove it shares the limit distribution of the previous
SDE. Let
\[
\check{\zeta}_t^{i,n} = \bigl( \check{P}^{i,n}_1(t)+
\check{P}^{i,n}_2(t) \bigr) + \int_0^t
\bigl\langle\check{Q}^{i;n}_1(s), \check{
\zeta}_s^{n} \bigr\rangle \,ds + \sum
_{l=1}^{d_x}\int_0^t
\bigl\langle\check{Q}^{i,l;n}_2(s), \check{\zeta
}_s^{n} \bigr\rangle \,dW^l_s,
\]
where
\begin{eqnarray*}
\check{Q}^{i;n}_1(s) &=&  \nabla F^{i}(X_s);\qquad
\check{Q}^{i,l;n}_2(s) = \nabla G^{i,l}(X_s);
\\
\check{P}^{i;n}_1(s)& = & \frac{1}{2}\int
_0^t n^r \bigl\langle\nabla
F^i (X_s), F(X_s) \bigr\rangle
\,dA^{0;n}
\\
&&{} + \sum_{l=1}^{d_x}\int
_0^t n^{r} \bigl\langle\nabla
F^i (X_s ), G^{\cdot
,l}(X_s) \bigr
\rangle \,dA_s^{1;l;n}
\\
&&{}+ \int_0^t n^{r} \bigl\langle
\nabla G^{i,l} (X_s), F(X_s) \bigr\rangle
\,dB_s^{0;l,n}
\\
&&{}+ \sum_{k,l=1}^{d_x} \frac{1}{\sqrt{2}} \int
_0^t n^{r} \bigl\langle\nabla
G^{i,l} (X_s) G^{\cdot,k}(X_s) \bigr
\rangle \,dB_{s}^{2;k,l,n},
\\
\check{P}^{i;n}_2(s) &= & n^{r-\beta}\int
_0^t \sum_{j,k,l=1}^{d_x}
\varphi_{G}^{i,j,l,k}(X_s) \,dB_s^{3;l,k,j,n}+n^{r-\beta}
\int_0^t \sum_{j=1}^{d_x}
R_{G}^{i,j}(X_s) \,dW^j_s
\\
&&{} + n^{r-\beta}\int_0^t \sum
_{j=1}^{d_x}\varphi_{F}^{i,j}(X_s)
\,dB_s^{1;j,n} + n^{r-\beta} \int_0^t
R_{F}^{i}(X_s) \,ds,
\end{eqnarray*}
for $R_F,R_G,\varphi_F,\varphi_G$ defined in Proposition~\ref
{PropMainDecreasingStepPqCLT}. By \ref{HypSlowScale}, $F,G$ are
bounded; by Lemma~\ref{LemPropertiesAverage}, $\nabla F$ and $\nabla G$ are well defined
and bounded and have bounded derivatives; and from the definition of
$R_F,R_G,\varphi_F,\varphi_G$ are $C_b^1$.

Note that \eqref{EqTupleConv3} in Proposition~\ref
{PropConvergenceOfTuple} gives us goodness and convergence of the
tuple $(n^r A^{0,n},n^r A^{1,n}$, $n^r B^{0,n},B^{1,n},n^r B^{2,n},
B^{3,n})$. Hence, by virtue of Theorem~5.4 in \citet{kurtzweak1991}
$\check{\zeta}^{n}(\cdot\wedge\tau_K^{n})$ is tight and any limit
point will satisfy \eqref{EqlimiDistSDE} on the interval $[0,\tau_K]$
where $\tau_K = (\inf\{t\dvtx|\zeta(t)|>K\} \wedge T)$. Moreover
\[
\sup_{0\leq s \leq\tau_K } \bigl\|\check{P}^{n}_s\bigr\|,\qquad \sup
_{0\leq s \leq
\tau_K }\bigl \|\check{Q}_1^{n}(s)\bigr\|, \qquad \sup
_{0\leq s \leq\tau_K } \bigl\| \check{Q}_2^{n}(s)\bigr\|
\]
are tight.

(c) We prove now that both $\zeta^n$ and $\check{\zeta}^n$ have the
same limit on the interval $[0,\tau_K]$. By Theorem~\ref
{ThmConvergenceChangeCoefficients}, it suffices to prove that sup norm
of the difference of the coefficients converge in probability. By
Theorem~\ref{TheoGeneralConvergence} the regularity properties of $F$
and the mean value theorem we have
\[
{{\mathbb{E}}} \Bigl[{\sup_{0\leq t \leq\tau_K} \bigl| Q_1^{i,n}(t)
- \check {Q}_1^{i,n}(t)\bigr|}\Bigr] \leq{{\mathbb{E}}} \Bigl[{
\sup_{x} \bigl| D^2 F(x) \bigr| \sup_{0\leq t \leq\tau_K}
\bigl| X_t - \ApproxVar{X}_t^n\bigr|}\Bigr]
\rightarrow0.
\]

The terms of $Q_2^{n},P_1^{n}$ are treated in the same way. On the
other hand, we get from Corollary~\ref{CorControlPowerAndDifference},
Proposition~\ref{PropAPriori}, and Burckholder--Davis--Gundy
inequality that
\begin{eqnarray*}
&&\sup_{0\leq t\leq T}\Biggl| n^{r-\beta}\int_0^t
\Biggl[ n^{\beta} \bigl( F^i(\cdot) - \DecrAvg{F}^{i;n}
(\cdot, \tbar[s] ) \bigr) - \sum_{j=1}^{d_x}
\varphi_F^{i,j}(\cdot) \nu^{j;n}_{\tbar[s]}
- R_{F}^i(\cdot) \Biggr] \bigl(\ApproxVar{X}^{n}_{\tbar[s]}
\bigr) \,ds\Biggr|,
\\
&&\sup_{0\leq t\leq T} \Biggl| n^{r-\beta}\int_0^t
\Biggl[ n^{\beta} \bigl( G^{i,j}(\cdot) - \DecrAvg{G}^{i,j;n}
(\cdot,\tbar[s] ) \bigr) \\
&&\hspace*{71pt}{}- \sum_{j,k=1}^{d_x}
\varphi_G^{i,j,l,k}(\cdot) \nu^{l,k;n}_{\tbar[s]}
- R_{G}^{i,j}(\cdot) \Biggr] \bigl(\ApproxVar{X}^{n}_{\tbar[s]}
\bigr) \,dW^l_s \Biggr|
\end{eqnarray*}
are tight and converge to zero.

Thus, by Theorem~\ref{ThmConvergenceChangeCoefficients} we will have
that $\zeta^{i;n}$ and $\check{\zeta}^{i;n}$ will converge to the
same limit.

(d) Finally, note that $\tau^{n}_K\rightarrow\infty$ and $\tau_K
\rightarrow\infty$, proving our claim in the full SDE case.

(e) To prove (i) it suffices to
follow the same approach. We obtain an equivalent development for the
ODE with random coefficients case (replacing by zero all the
``g-terms''). The rest of the proof proceeds as before, this time using
\eqref{EqTupleConv1} for the weak convergence of the tuple.
\end{pf*}

\section{The EMsDS algorithm}
\label{SecEMSDSalgorithm}

Given the error expansion for the decreasing step algorithm presented
in Proposition~\ref{PropErrorExpansion}, it seems natural to explore
if a Richardson--Romberg extrapolation may be used to obtain the
approximation with the same convergence properties we have proven. The
idea of such a procedure is to decrease the complexity by performing a
linear combination of two (or more) realizations of the algorithm with
carefully chosen parameters. We borrow here the procedure as defined in
\citet{lemaireestimation2005}.

Let $\lambda$ be a positive real. If $\{\gamma_k\}$ is a sequence of
steps satisfying \ref{HypSteps}, the sequence $\gamma
_k^\lambda
:= \frac{\gamma_k}{\lambda}$ will also satisfy \ref{HypSteps}.
We will denote $\Gamma^\lambda_M$ and $\Gamma^{\lambda,[r]}_M$ the
sum of the $\gamma_k^\lambda$ and its power as before.

Let us denote by $\DecrAvg{F}^{\lambda,M}(x,q)$ the approximation as
defined in \eqref{EqAverageEstimator} when the coefficients $\{\gamma
_k^\lambda\}_{k\in\NN^*}$ are used.

With $\varsigma$ given as in \eqref{EqDefinitionLbar}, let us define
the \emph{extrapolated approximation estimator} as
%
\begin{equation}
\label{EqDefinitionFHat} \DecrAvgExtrap{F}^{\lambda;M(n)}(x,q) = \frac{1}{\lambda^{\varsigma
/2-1}-1}
\bigl(\lambda^{\varsigma/2-1} \DecrAvg{F}^{\lambda, M(n)
}(x,q) -
\DecrAvg{F}^{ M(n) }(x,q) \bigr).
\end{equation}

The first question we might ask is if estimator \eqref
{EqDefinitionFHat} does converge to the actual ergodic average, and
what type of properties it inherits. To clarify the situation consider
an extension of \eqref{EqErgodicDiffusion}. Let $\vec
{Y}^x=(Y^{1;x},Y^{2;x})^*$ with
%
\begin{eqnarray}
\label{EqExtendedSystem} %
Y_{t}^{1;x} &=&
y_0^1 + \int_0^t
\frac{b(x,Y_s^{1;x})}{\lambda} \,ds + \int_0^t
\frac{\sigma(x, Y_s^{1;x})}{\sqrt{\lambda}} \,d\ApproxVarExtrap{W}^1_s,
\nonumber
\\[-8pt]
\\[-8pt]
\nonumber
Y_t^{2;x}& =& y_0^2 + \int
_0^t b \bigl(x,Y_s^{2;x}
\bigr) \,ds + \int_0^t \sigma \bigl(x,
Y_s^{2;x} \bigr) \,d\ApproxVarExtrap{W}^2_s.
\end{eqnarray}
If $\ApproxVarExtrap{W}^1$ and $\ApproxVarExtrap{W}^2$ are independent,
then this system satisfies \ref{HypFastScale} with a unique
invariant measure
defined by $\vec{\mu}^x(d\vec{y}) = \mu^x(dy^1)\mu^x(dy^2)$. If we define
%
\begin{equation}
\label{EqExtendedF} \vec{f}(x,\vec{y}):= \frac{1}{\lambda
^{\varsigma/2-1}-1} \bigl(\lambda
^{\varsigma/2-1} f \bigl(x,y^1 \bigr)-f \bigl(x,y^2
\bigr) \bigr),
\end{equation}
and defining in an analogous way $\vec{h}$, then it can be seen that
$\vec{f}, \vec{g}, \vec{h}:=\vec{g}\vec{g}^*$ satisfy
\ref{HypSlowScale}. Moreover
if we apply the decreasing step algorithm to $\vec{f}$ (resp., $\vec
{h}$) in the extended framework, we obtain the expression \eqref
{EqDefinitionFHat}. \emph{Hence, we conclude that the EMsDS algorithm
is equivalent to the MsDS algorithm applied to an extended system.}

Let us denote by $\ApproxVarExtrap{X}^{n}$ the approximation of the
diffusion $X$ using the extrapolated version of the algorithm. In view
of the discussion we presented before, the following result is mainly a
corollary of Theorems \ref{TheoGeneralConvergence} and \ref
{TheoAbstractLimitDistribution}, and extends the main Theorem to the
extrapolation algorithm. It shows the advantage of using the EMsDS
algorithm: assuming higher regularity, all the properties of the MsDS
algorithm are conserved, but the extrapolated version \emph{allows a
lower value for $\theta$ in the definition of the sequence $\gamma_k =
\gamma_0 k^{-\theta}$}. More precisely we pass from $1/2$ to $1/3$ in
the ODE case and from $1/3$ to $1/5$ in the SDE case as minimal $\theta
$ values. As a consequence of this reduction, the complexity of the
modified version is in general asymptotically lower than that of the
nonextrapolated version (refer to the efficiency analysis on
Section~\ref{SubsecEfficiency}).

\begin{MyTheorem}
\label{TheoShortMainDecreasingStepInterp}
Let $0<\theta<1$, $\gamma_1 \in\RR^+$ and $\gamma_k = \gamma_1
k^{-\theta}$. Assume \ref{HypFastScale} and
\ref{HypSlowScale}, $M(n)$ defined as in Theorem~\ref{TheoShortMainDecreasingStep}, and assume in addition that $r_y>
5$. Let $\ApproxVarExtrap{X}^n$ be the approximated diffusion where we
replace the ergodic estimator \eqref{EqAverageEstimator} by~\eqref
{EqDefinitionFHat}.
\begin{enumerate}[(i)]
\item[(i)] (Strong convergence). There exists a constant $K$ such that
\begin{itemize}
\item Case $g\equiv0$ (ODE with random coefficients):
\[
{{\mathbb{E}}} \Bigl[{\sup_{0\leq t \leq T}\bigl|X_t -
\ApproxVarExtrap {X}_t^{n} \bigr|^2 }\Bigr] \leq K
n^{-2[(1-\theta)\wedge2\theta]/(1-\theta)}.
\]
\item(Full SDE case):
\[
{{\mathbb{E}}} \Bigl[{\sup_{0\leq t \leq T}\bigl|X_t -
\ApproxVarExtrap {X}_t^{n}\bigr |^2 }\Bigr] \leq K
n^{-[(1-\theta)\wedge4\theta]/(1-\theta)}.
\]
\end{itemize}
\item[(ii)] (Limit distribution). Assume in addition that $r^y \geq8$,
and define
\begin{eqnarray*}
\ApproxVarExtrap{C}_\varphi&:= &\frac{(\lambda^{3}+1)^{1/2}}{\lambda
-1};\qquad \ApproxVarExtrap{C}_1:=
\frac{\gamma_0^2 (1-\theta)
M_1^{-2\theta
} }{1-3\theta};
\\
\ApproxVarExtrap{\varphi}_F(x) &:=& \mathbf{1}_{\{\theta= 1/5\}}
\ApproxVarExtrap{C}_\varphi\sqrt{\ApproxVarExtrap{\Phi}_F(x)};\qquad
\ApproxVarExtrap{\varphi}_G(x):= \mathbf{1}_{\{\theta= 1/5\}}
\ApproxVarExtrap{C}_\varphi\sqrt{\ApproxVarExtrap{\Phi}_G(x)};
\\
\ApproxVarExtrap{R}_{F}^{i}(x) &:=& \mathbf{1}_{\{\theta\geq1/5\}}
\ApproxVarExtrap{C}_1 \bigl(1-\lambda^{-1} \bigr) \int
\bar{v}_{F^i}^{\varsigma
+2,r^y}(x,y) \mu^x(dy);
\\
\ApproxVarExtrap{R}_{H}^{i,j}(x) &:=& \mathbf{1}_{\{\theta\geq1/5\}}
\ApproxVarExtrap{C}_1 \bigl(1-\lambda^{-1} \bigr)\int\bar
{v}_{H^{i,j}}^{\varsigma
+2,r^y}(x,y) \mu^x(dy).
\end{eqnarray*}
\begin{itemize}
\item{[ODE case: $G(x)\equiv0$].}
If $\theta\geq1/3 $, then $\ApproxVarExtrap{\zeta}^{n}:= n (X_t
- \ApproxVarExtrap{X}^{n} ) $ satisfies the limit distribution
result given in Theorem~\ref{TheoShortMainDecreasingStep}\textup{(a)} with
new coefficients $ \ApproxVarExtrap
{\varphi}^F$ instead of $\varphi^F$.
\item
(SDE case). If $\theta\geq
1/5$, then $\ApproxVarExtrap{\zeta}^{n}:= n^{1/2} (X_t -
\ApproxVarExtrap{X}^{n} ) $ satisfies the limit distribution
result given in Theorem~\ref{TheoShortMainDecreasingStep}\textup{(b)} with the coefficients $\ApproxVarExtrap
{R}^F,\ApproxVarExtrap{R}^G,\ApproxVarExtrap{\varphi}^F$ and
$\ApproxVarExtrap{\varphi}^G$ instead of $R^F$, $R^G$, $\varphi_F$ and
$\varphi_g$, respectively.
\end{itemize}
\end{enumerate}
\end{MyTheorem}

\begin{pf*}{Proof of Theorem~\ref{TheoShortMainDecreasingStepInterp}}
We will deduce the proof only for the full SDE case the other case
being analogous. We assume that $\varsigma=4$, which is the most
common case.

\begin{longlist}[(a)]
\item[(a)] As in the proof of Theorem~\ref{TheoShortMainDecreasingStep},
the sequence of coefficients satisfies \ref{HypSteps}.
Moreover, the EMsDS
algorithm is the MsDS algorithm applied to an extended system, and
hence the strong convergence and limit distribution properties are a
consequence from Theorems \ref{TheoGeneralConvergence} and \ref
{TheoAbstractLimitDistribution}: it remains just to express the values
of the functions and constants appearing in Propositions \ref
{PropControlDecreasingStep} and \ref{PropMainDecreasingStepPqCLT} in
terms of the original system.

Indeed, recall that
%
\begin{eqnarray}\qquad
\label{EqDefinitionDerivativeHat} \vec{b}(x,\vec{y}) = %
\pmatrix{ \lambda^{-1}b
\bigl(x,y^1 \bigr) \vspace*{2pt}
\cr
b \bigl(x,y^2
\bigr) } %
;\qquad \vec{\sigma}(x,\vec{y}) = %
\pmatrix{
\lambda^{-1/2} \sigma \bigl(x,y^1 \bigr) & 0 \vspace*{2pt}
\cr
0 & \sigma \bigl(x,y^2 \bigr) }.
\end{eqnarray}
By (i) in Proposition~\ref{PropControlDecreasingStep} applied to the
extended problem [i.e., for the system~\eqref{EqExtendedSystem} and
$\vec{f}$ defined in \eqref{EqExtendedF}], we have a solution for the
extended centered Poisson equation given by
\[
\vec{\phi}_{F^i}(x,\vec{y}) = (\lambda-1)^{-1} \bigl(
\lambda^{2}\phi_{F^i} \bigl(x,y^1 \bigr) -
\phi_{F^i} \bigl(x,y^2 \bigr) \bigr),
\]
that is, the solution of equation \eqref{EqCenteredPoisson} with
function $F^i$ under the extended set-up is a linear combination of the
solution in the original set-up. Thus, for any $j>0$,
%
\begin{equation}
\label{EqDerivativePhiHat} D^j_y \vec{\phi}_{F^i}(x,
\vec{y}) = \frac{1}{\lambda-1} %
\pmatrix{ \lambda^{2}
D^j_y\phi_{F^i} \bigl(x,y^1 \bigr)
\vspace*{2pt}
\cr
- D^j_y \phi_{F^i}
\bigl(x,y^2 \bigr) } %
.
\end{equation}
It follows that
\begin{eqnarray*}
&& D^j_y\vec{\phi}_{F^i}(x,\vec{y}){{
\mathbb{E}}} \bigl[{ \bigl\langle \vec{b}(x,\vec{y})^{\otimes(l-j)}, \bigl(\vec{
\sigma}(x,\vec{y}) U_1^0 \bigr)^{\otimes(2j-l)} \bigr
\rangle}\bigr]
\\
&&\qquad = \frac{\lambda^{2}}{\lambda-1} D^j_y\phi_{F^i}
\bigl(x,y^1 \bigr){{\mathbb{E}}} \biggl[{ \biggl\langle \biggl(
\frac
{b(x,y^1)}{\lambda} \biggr)^{\otimes(l-j)}, \biggl(\frac {\sigma
(x,y^1)}{\sqrt{\lambda}}
U_1^0 \biggr)^{\otimes(2j-l)} \biggr\rangle}\biggr]
\\
& &\qquad\quad{}-\frac{1}{\lambda-1} D^j_y\phi_{F^i}
\bigl(x,y^2 \bigr){{\mathbb{E}}} \bigl[{ \bigl\langle\vec{b}
\bigl(x,y^2 \bigr)^{\otimes(l-j)}, \bigl(\sigma \bigl(x,y^2
\bigr) U_1^0 \bigr)^{\otimes(2j-l)} \bigr\rangle}
\bigr].
\end{eqnarray*}
Therefore
%
\begin{equation}
\label{Eqvarrowbar} \vec{\bar{v}}^{l,r^y}_{F^i} = \biggl(
\frac{\lambda^{(4-l)/2}-1}{\lambda
-1} \biggr)\bar{v}^{l,r^y}_{F^i},
\end{equation}
and we deduce that the terms of the error expansion will be zero for
$l\leq5$.

\item[(b)] Let $\vec{\varsigma}$ be defined by \eqref{EqDefinitionLbar}
under the extended setup. From \eqref{Eqvarrowbar} we conclude that
$\vec{\varsigma}\geq6$, being $\vec{\varsigma}=6$ the worst case.
Hence, we deduce that defining
\[
\beta_0 = \frac{1}{2},\qquad \ApproxVarExtrap{\beta_1}
=\frac{2\theta
}{1-\theta},\qquad \hat{\beta}_2 =\frac{\theta}{1-\theta},
\]
then
\begin{eqnarray*}
\Gamma_M &\approx&\frac{\gamma_0 M_1^{1-\theta}n}{1-\theta},\qquad \frac
{\Gamma_M^{[3]}}{\Gamma_M} \approx
\frac{\gamma_0^2 (1-\theta)
M_1^{-2\theta} n^{-{2\theta}/{(1-\theta)}}}{1-3\theta},
\\
\frac{\Gamma_M^{[4]}}{\Gamma_M^{[3]}} &\approx&\frac{\gamma_0
(1-3\theta) M_1^{-\theta} n^{-{\theta}/{(1-\theta)}}}{1-4\theta},
\end{eqnarray*}
and so, $\beta_1, \ApproxVarExtrap{\beta}_2,\ApproxVarExtrap{\beta}_3$
are the coefficients appearing in Proposition~\ref
{PropControlDecreasingStep} applied to this setup. We conclude as well
that $\ApproxVarExtrap{R}_{F}^{i}$ is the function appearing in
Proposition~\ref{PropMainDecreasingStepPqCLT}. Similar developments
for $H$ allow us to extend the conclusion to~$\ApproxVarExtrap{R}_{H}^{i,j}$.

\item[(c)] Finally, looking at the definition of $\varphi_F$ and $\Phi_F$
from Proposition~\ref{PropMainDecreasingStepPqCLT} and~\eqref
{EqDerivativePhiHat} we get that
\begin{eqnarray*}
&&\ApproxVarExtrap{\Phi}_F^{i,j}(x)=\frac{C_0^{-1}}{(\lambda-1)^2}
\biggl( \lambda^{2} \int \bigl\langle\sigma^* D_y
\phi_{F^i}, \sigma^* D_y \phi_{F^j} \bigr\rangle
\bigl(x,y^1 \bigr) \mu^x \bigl(dy^1 \bigr)
\\
&&\hspace*{93pt}\qquad{}+\int \bigl\langle\sigma^* D_y \phi_{F^i}, \sigma^*
D_y \phi_{F^j} \bigr\rangle \bigl(x,y^2 \bigr)
\mu^x \bigl(dy^2 \bigr) \biggr);
\end{eqnarray*}
that is, $\ApproxVarExtrap{\Phi}_F(x) = (\lambda^{2}+1)(\lambda-1)^{-2}
\Phi_F(x)$. We get a similar result for $\ApproxVarExtrap{\Phi}_G$. We
obtain the value $\ApproxVarExtrap{C}_\varphi$ given in the statement.
The claim follows.\quad\qed
\end{longlist}
\noqed\end{pf*}

\begin{MyRemark}
$\ApproxVarExtrap{C}_\varphi$ is a constant multiplying the
uncertainty coming from the decreasing step estimator. Since we would
like this quantity as small as possible, having an explicit value for
$\ApproxVarExtrap{C}_\varphi$ is very useful from a numerical point of
view: we can choose $\lambda$ to minimize $\ApproxVarExtrap
{C}_\varphi
$. We get
\[
\lambda_* = 1+(\sqrt{3}+1)^{1/3}+(\sqrt{3}+1)^{-1/3}
\approx3.196
\]
inducing
$\ApproxVarExtrap{C}_\varphi\approx2.64 $. This is the initial
additional cost that has to be paid for the extrapolation, making the
EMsDS algorithm useful for large $n$, where the reduction in complexity
of the EMsDS is enough to compensate for the higher error.
\end{MyRemark}

\section{Numerical results}
\label{Secnumresults}
\subsection{Efficiency analysis}
\label{SubsecEfficiency}
We can approximate the execution time of both algorithms, the original
and extrapolated versions of the algorithm, by estimating the total
number of operations needed to perform one path\vadjust{\goodbreak} approximation of the
effective equation \eqref{EqEffectiveEquation}. Note that since both
algorithms share the same structure, a~similar analysis is valid for
both of them: the total cost $\kappa(n)$ of the algorithm with $n$
steps may be written as
\[
\kappa(n)= \bigl[\kappa_1(n,d_x,d_y)+
\kappa_2(d_x) \bigr]n,
\]
where $\kappa_1$ stands for the cost coefficient estimation at each
step of the decreasing Euler, and $\kappa_2$ for the cost of
calculating the Euler iteration. The latter will be of order $O(d_x)$
in the ODE case and $O(d_x^2)$ for the SDE case.

Let us focus now on $\kappa_1$. Both algorithms perform
$M_1n^{1/(1-\theta)}$ iterations for approximating the diffusion
$\DecrAvg{Y}$ and the calculation of estimators $\DecrAvg{F},
\DecrAvg
{G}$. For the MsDS algorithm, each one of these iterations has a cost
of $O(d_y d_x)$ in the $\mathit{ODE}$ case, or $O(d_y d_x^2)$ in the SDE case.
In the latter, we need also to perform a Cholesky decomposition with a
cost of $O(d_x^3)$ operations. Hence
\[
\kappa_1^{\mathrm{MsDS}} (n,d_x,d_y ) =
\cases{ O \bigl( d_yd_x n^{1/(1-\theta)}
\bigr), &\quad$\mbox{ in ODE case},$ \vspace*{2pt}
\cr
O \bigl( \bigl[d_yd_x^2
+ d_x^3 \bigr] n^{1/(1-\theta)} \bigr), &\quad$\mbox{ in SDE
case}.$}
\]
On the other hand, from the definition of the EMsDs algorithm, we get
$\kappa_1^{\mathrm{EMsDS}} \leq\lambda\kappa_1^{\mathrm{MsDS}}$, and thus both share
the same order of complexity, with the only difference that \emph
{$\theta$ is allowed to be smaller in the extrapolated algorithm}.

It may be more interesting to compare the \emph{efficiency} of both
algorithms, that is, the time spent to obtain a given error tolerance
$\Delta$. We have from Theorems \ref{TheoShortMainDecreasingStep} and~\ref{TheoShortMainDecreasingStepInterp} that $\Delta(n):= O( n^{-1} )
$ for the ODE, and $\Delta(n):= O( n^{-1/2} )$ for the SDE case.
Replacing the minimum possible $\theta$ values we obtain the complexity
figures given in Table~\ref{TabEfficiency}.
%
\begin{table}[b]
\caption{Minimal efficiency (operations for fixed error) of the basic
and extrapolated algorithm for ODE and full SDE cases}\label{TabEfficiency}
\begin{tabular*}{\textwidth}{@{\extracolsep{\fill}}lcccc@{}}
\hline
& \textbf{ODE} & \textbf{ODE (extrapol.)} & \textbf{SDE} & \textbf{SDE (extrapol.)}\\
\hline
$\theta_{\min}$ & $1/2$ & $1/3$ & $1/3$ & $1/5$\\
$\tau_{\min}(\Delta)$ & $O(d_yd_x \Delta^{-3}) $ & $O(d_xd_y \Delta
^{-2.5})$ & $O([d_x^2d_y+d_x^3] \Delta^{-5})$ & $O([d_x^2d_y+d_x^3]
\Delta^{-4.5})$\\
\hline
\end{tabular*}
\end{table}

How do these figures compare with a straightforward Euler scheme
applied to the original system? For the ODE case, an Euler scheme
implemented for the original system \eqref{EqTheSystem} would require
a total of $(dx+dy)\varepsilon^{-1} \Delta^{-2}$ operations. Then the MsDS
algorithm is more efficient if $\varepsilon< \Delta(d_x\vee d_y)^{-1}$,
and the EMsDS if $\varepsilon< \Delta^{1/2} (d_x\vee d_y)^{-1}$. With
respect to the algorithm presented in \citet{eanalysis2005}, the
efficiency is equivalent to the one obtained when using a weak scheme
of order one for approximating the ergodic averages. The advantage of
our method is that we have in addition to the rate of convergence an
expression for a C.L.T. type result.

In the SDE case, on the other hand, the proposed algorithm will be
advantageous in the case in which $\varepsilon< \Delta^3 (d_x\vee
d_y)^{-1}$ for the MsDS version, and $\varepsilon< \Delta
^{2.5}(d_x\vee
d_y)^{-1}$ for the EMsDS. In other words, our proposed algorithms will
be more efficient in our regime of interest of a strong scale
separation (i.e., when $\varepsilon\rightarrow0$). It should be remarked
that the SDE case is not explicitly studied for the algorithm in \citet
{eanalysis2005}.

\subsection{Numerical tests}
\label{SubSecNumericalTests}
\subsubsection{A toy problem}
Let us illustrate the main features of the algorithm by evaluating its
behavior when used for solving a toy system for which we are able to
obtain an exact solution. Consider
\[
dY^x_t = \bigl( \bigl(|x|^2+1
\bigr)^{-1/2} - Y^x_t \bigr) + \sqrt{2}\,d
\tilde{W}_t,
\]
which is an Ornstein--Uhlenbeck system having a unique invariant
measure with normal distribution with mean $(|x|^2+1)^{-1/2}$ and
variance $1$, and define the SDE system
\[
dX_t = F(X_t) \,dt + G(X_t) \,dW_t,
\]
with
\begin{eqnarray*}
f(x,y)& :=& %
\pmatrix{ 1+y- \bigl(|x|^2+1
\bigr)^{-1/2} \vspace*{2pt}
\cr
1 } %
; \\
 g(x,y) &:=& \sqrt{
\frac{|x|^2+1}{2|x|^2+3} \bigl(y^2+1 \bigr)} %
\pmatrix{ 1 & 0
\vspace*{2pt}
\cr
1 & 1} %
,
\end{eqnarray*}
with $F,G$ defined as before and where $\tilde{W}$ is a real Brownian
motion independent of the planar Brownian motion $W$.
The form of the assumed coefficients is chosen to satisfy the
regularity and uniform bound hypothesis in \ref{HypSlowScale}
and \ref{HypFastScale} and to
give a simple effective equation expression. In fact, it is easily
verified that the exact effective equation is
\[
X_s = %
\pmatrix{ x_0^1 + s +
W^1_s\vspace*{2pt}
\cr
x_0^2
+ s + W^1_s+W^2_s }
.
\]

\begin{figure}

\includegraphics{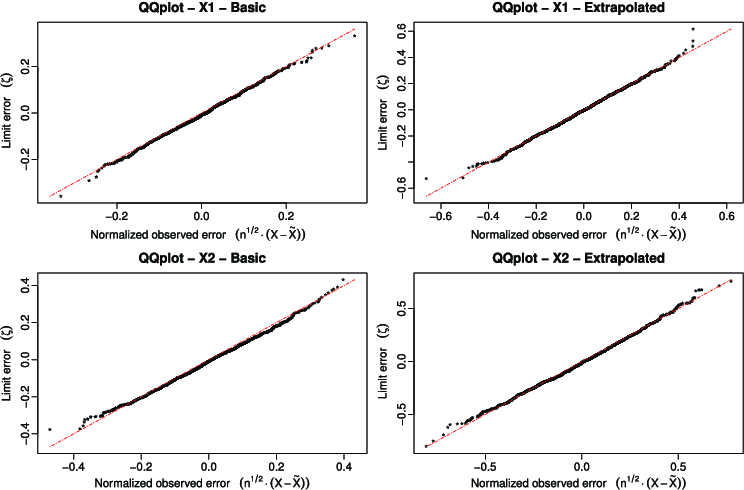}

\caption{Q--Q plot comparing the rescaled errors in the simulation with
$n=510$ and the theoretical limit distribution (the reference line
represents a perfect match). Left: SDE decreasing step. Right: SDE
interpolated.}
\label{figqqplot}
\end{figure}

\begin{figure}

\includegraphics{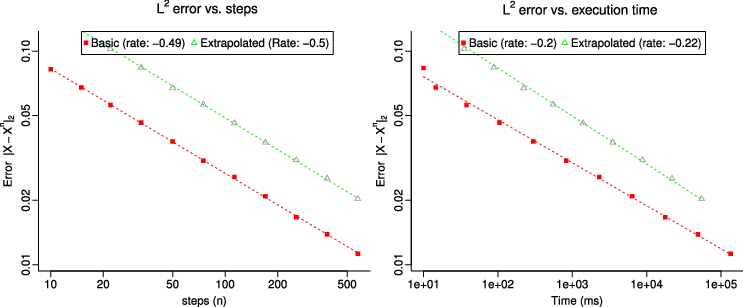}

\caption{Left: $L_2$ error as a function of steps for the SDE case
(log--log scale). Note that the estimated values for the slopes verify
the rate of convergence for the algorithm in both implementations.
Right: $L_2$ error as a function of execution time for the SDE case
(log--log scale). Although a higher price must be payed for a small step
number, the slope difference signals a change in the asymptotic order
of convergence.}
\label{figStrongErrorSteps}
\end{figure}

We will look at the numerical results of applying the decreasing step
with sequence $\gamma_k = k^{-1/3}$ and the EMsDS version with sequence
$\gamma_k = k^{-1/5}$ and $\lambda=3$. Let us examine the distribution
of the error at a fixed time $T=1$ (i.e., $\zeta= \ApproxVar{X}_1 -
X_1 $). Figure~\ref{figqqplot} shows a Q--Q plot of the rescaled
simulated errors $\sqrt{n} \zeta$ and the limit distribution error in
the studied cases. As shown, the empirical distributions obtained after
1600 simulations with $n = 510$ verify the expected limit behavior.

Figure~\ref{figStrongErrorSteps}. Left plots in a log--log scale the
evolution of the $L_2$ error
\[
\zeta_{L_2} = \Bigl({{\mathbb{E}}} \Bigl[{\sup_{0\leq t \leq T}
| \ApproxVar{X}_t - X_t |^2}\Bigr]
\Bigr)^{1/2}\vadjust{\goodbreak}
\]
in function of the number of steps $n$, comparing both versions of the
algorithm. The empirically obtained slope (close to $-0.5$ in both
cases) represents the power of the approximation and is the one
expected from the convergence theorems.

We show as well in Figure~\ref{figStrongErrorSteps} (right) a comparison
in the efficiency of both methods (measured as the error in terms of
the execution time) of each one of the algorithms. The effect of the
extrapolation in the cost of the algorithm is evidenced in the
difference in slope of the empirical plot for both algorithms. Note
that solving for $\Delta$ in Table~\ref{TabEfficiency} we get $\Delta
_{\mathrm{MsDS}} = O(\tau^{-0.2})$ and $\Delta_{\mathrm{EMsDS}} \approx O(\tau
^{-0.222})$, values that are retrieved in the numerical experiment. It
is worth observing the difference in the intercept of both lines,
showing that the higher slope comes with a cost in the initial error.
The conclusion drawn from the toy example may well be generalized: the
user should consider implementing the extrapolated version only when
requiring a very high precision on the approximation results.
\subsubsection{Pricing in finance} We apply now the algorithm to a
pricing problem in finance. Consider the mean-reverting corrected
Heston's stochastic volatility model presented in \citet
{fouquefast2011} and given by
\begin{eqnarray*}
dX_t &=& rX_t \,dt + \Sigma_t X_t
\,dW^{x}_t,
\\
 dY_t& =& \varepsilon^{-1} Z_t
(m-Y_t) \,dt + \nu\sqrt{2 Z_t \varepsilon^{-1} }
\,dW^y_t,
\\
 dZ_t &=& \kappa(\theta-Z_t)\,dt + \sigma
\sqrt{Z_t} \,dW^z_t,
\\
 \Sigma_t &=& \sqrt{Z_t} \bigl(1+Y_t^2
\bigr),
\end{eqnarray*}
where we assume $W_t^x,W^y_t,W^z_t$ are one-dimensional Brownian
motions with correlations $\rho_{xy}$, $\rho_{xz}$ and $\rho_{yz}$. We
suppose the model is already written in terms of the risk neutral
probability measure with known parameters and initial conditions given
in Table \ref{Tableparameters}.
%
\begin{table}[b]
\caption{Initial condition and parameters of the model}
\label{Tableparameters}
\begin{tabular*}{\textwidth}{@{\extracolsep{\fill}}lccccccccccc@{}}
\hline
$\bolds{x_0}$ & $\bolds{z_0}$ & $\bolds{y_0}$
& $\bolds{m}$& $\bolds{\nu}$& $\bolds{\kappa}$ &$\bolds{r}$&
$\bolds{\theta}$& $\bolds{\sigma}$
&$\bolds{\rho_{xy}}$ &$\bolds{\rho_{yz}}$ &$\bolds{\rho_{xz}}$ \\
\hline
100 &0.24 &0.06 & 0.06 & 1.0 & 1.0 &0.05 & 1.0 &0.39 & 0 & 0 & $-0.33$\\
\hline
\end{tabular*}\vspace*{-3pt}
\end{table}
We are interested in pricing several types of options depending on the
whole trajectory on this model. For this test, we price a floating
strike Asian call (the payoff being
$AC_{\mathrm{float}}=S_T-T^{-1} \int S_t \,dt$) and a lookback call with floating strike (with payoff
$LC_{\mathrm{float}}=S_T-S_{\min}$).

In this test, we compare the algorithm with a simple Euler scheme with
different values for $\varepsilon$. We carry out 6000 Monte Carlo
simulations. The results are presented in Table~\ref{Tablesimulation}.

\begin{table}[]
\caption{Simulation values}
\label{Tablesimulation}
\begin{tabular*}{\textwidth}{@{\extracolsep{\fill}}lcccccc@{}}
\hline
\textbf{Method} & $\bolds{\varepsilon}$ & $\bolds{n}$ & $\bolds{M(n)}$
& $\bolds{n \times M(n)}$& \textbf{Asian} & \textbf{Lookback}
\\
\hline
Euler & $10^{-3}$ & $5 \times10^6$ & $1$ & $5 \times10^4$ & 40,988 &
81,591\\
Euler & $10^{-3}$ & $ 10^7$ & $1$ & $1 \times10^5$ & 40,503 & 81,256\\
Euler & $10^{-3}$ & $2\times10^7$ & $1$ & $2 \times10^5$ & 40,086 &
80,769\\
Euler & $10^{-4}$ & $5 \times10^6$ & $1$ & $5 \times10^5$ & 22,091 &
54,119\\
Euler & $10^{-4}$ & $ 10^7$ & $1$ & $1 \times10^6$ & 21,897 & 53,806\\
Euler & $10^{-4}$ & $2\times10^7$ & $1$ & $2 \times10^6$ & 20,908 &
52,095\\
Euler & $10^{-5}$ & $5 \times10^6$ & $1$ & $5 \times10^6$ & 18,203 &
45,947\\
Euler & $10^{-5}$ & $ 10^7$ & $1$ & $1 \times10^7$ & 15,164 & 39,123\\
Euler & $10^{-5}$ & $2\times10^7$ & $1$ & $2 \times10^7$ & 20,659 &
51,240\\[3pt]
MsDS & -- & 50 & 3540 & $1.77 \times10^5$ & 20,738 & 47,920\\
MsDS & -- & 100 & 10,010 & $1 \times10^6$ & 20,681 & 48,841\\
MsDS & -- & 200 & 28,290 & $5.66 \times10^6$ & 20,669 & 49,557\\
\hline
\end{tabular*}
\end{table}

Note that the system does not satisfy all the hypothesis
\ref{HypFastScale}
and
\ref{HypSlowScale}, particularly it fails to satisfy the
boundedness of the
coefficients with respect to the slow variables, and the uniform
ellipticity hypothesis. Nevertheless, the MsDS algorithm seems to work
even under these relaxed conditions, and, in addition, appears to be
more stable than the algorithm using small values of $\varepsilon$. Note
as well that for similar values of total operations [represented by the
column $n\times M(n)$], the MsDS algorithm gives better results.

\begin{appendix}\label{app}
\section{Technical results}
\label{SecTechnicalResults}

\subsection{Weak convergence of tuples}\label{SubSecWeakConvTuples}

\mbox{}

\begin{pf*}{Proof of Proposition~\ref{PropConvergenceOfTuple}}
(a) Let us start by proving \eqref{EqTupleConv1}. Note that the
approximations defined by \eqref{Eqcontinuousapproximated} are
defined in the same sample space of the effective equation \eqref
{EqEffectiveEquation} and that we have, thanks to
Theorem~\ref{TheoGeneralConvergence}, that
\[
\sup_{0\leq s\leq t}\bigl |\ApproxVar{X}^{n}_s -
X_s\bigr| \mathop{\stackrel{P} {\longrightarrow}}0.
\]
Hence,
%
\begin{equation}
\label{EqConvTriple} \bigl(X,\ApproxVar{X}^{n},W \bigr)\Rightarrow(X,X,W).
\end{equation}
Now, $nA^{0;n}$ is deterministic, continuous and
\[
\lim_{n\rightarrow\infty} n A^{0;n}_t = \lim
_{n\rightarrow\infty} 2n \biggl(\frac{t^2}{2} - \frac{\lfloor nt
\rfloor(\lfloor nt \rfloor-
1)}{2n^2}
\biggr) = t,
\]
and the convergence is uniform in $t$.

On the other hand, we can easily verify that for any $t$, and $1\leq i
\leq d_x$,
\[
\sqrt{n}B^{1;i,n}_t = \sqrt{n} \int_0^t
\nu^{i;n}_{\tbar[s]} \,ds = \sum_{i=0}^{\lfloor nt \rfloor}
\frac{1}{\sqrt{n}} \nu^{i;n}_{\tbar
[s]} + \frac{nt - \lfloor nt \rfloor}{n^{3/2}}
\nu^{i;n}_{\tbar},
\]
but by the Cauchy--Schwarz inequality we have
\[
{{\mathbb{E}}} \biggl[{\biggl|\sup_{0\leq t\leq T} \biggl(\frac{nt - \lfloor nt
\rfloor }{n^{3/2}}
\nu^{i;n}_{\tbar} \biggr) \biggr|^2}\biggr] \leq{{\mathbb
{E}}} \Biggl[{\Biggl|\sum_{k=1}^{n}
\frac{1}{n^{3/2}}\bigl |\nu^{i;n}_{t_k}\bigr| \Biggr|^2}\Biggr]
\leq{{\mathbb{E}}} \biggl[{\frac {n}{n^3}\sum{\bigl|\nu^{i;n}_{t_k}\bigr|^2}}
\biggr] \rightarrow0.
\]
Then it suffices to study the convergence of the Gaussian martingale\break
$\sum_{i=0}^{\lfloor nt \rfloor} n^{-1/2} \times  \nu^{i;n}_{\tbar[s]} $. Let
$0\leq j\leq d_x$. Then the independence properties and an application
of a multi-dimensional C.L.T. gives us that
\[
\Biggl\langle\sum_{i=0}^{\lfloor nt \rfloor}
\frac{1}{\sqrt{n}} \nu^{i;n}_{\tbar[s]},\sum
_{i=0}^{\lfloor nt \rfloor} \frac{1}{\sqrt{n}} \nu^{j;n}_{\tbar[s]}
\Biggr\rangle= \frac{1}{n} \sum_{i=0}^{\lfloor
nt \rfloor}
\nu^{i;n}_{\tbar[s]}\nu^{j;n}_{\tbar[s]} \mathop{
\stackrel {P} {\longrightarrow}} \delta_{i=j}.
\]
We conclude that $B^{1}$ is (up to a modification) a Brownian motion
independent from $W$ and $X$, by remarking its Gaussian nature with
independent increments property and covariance matrix as the one of the
standard Brownian. Thus \eqref{EqTupleConv1} follows. Note that we
have shown property (\ref{PropertyStar}) as well, and consequently goodness of the sequence.

(b) To prove \eqref{EqTupleConv3} note first that $\sqrt
{n}B^{2;i,j;n}$ is a continuous martingale. In view of the results in
\citet{jacodcontinuous1997}, we examine the component-wise quadratic
variation. By standard techniques we find
\begin{eqnarray*}
\bigl\langle\sqrt{n} B^{2;i,j;n}, \sqrt{n}B^{2;i',j;n} \bigr
\rangle_t &= &2n\int_0^t
\bigl(W^i_s-W^i_{\tbar[s]} \bigr)
\bigl(W^{i'}_s-W^{i'}_{\tbar[s]} \bigr)
\,ds
\\
& \mathop{\stackrel{P} {\longrightarrow}}&\mathbf{1}_{\{i=i'\}}t,
\end{eqnarray*}
and due to independence we find as well that, taking, $j\neq j'$,
\begin{eqnarray*}
\bigl\langle\sqrt{n} B^{2;i,j;n}, W^j \bigr
\rangle_t &\mathop{\stackrel{P} {\longrightarrow}}&0,
\\
\bigl\langle\sqrt{n}B^{2;i,j;n}, \sqrt{n}B^{2;i,j';n} \bigr
\rangle_t &=& 0 ;\qquad  \bigl\langle\sqrt{n}B^{2;k;n},
W^j \bigr\rangle_t = 0.
\end{eqnarray*}
By Theorem~4-1 in \citet{jacodcontinuous1997}, $B^{2;n}$ convergences
stably in law toward $B^{2}$ a standard Brownian Motion independent
from $W$; for the definition of this type of convergence see \citet
{aldousmixing1978} or \citet{jacodcontinuous1997}. Since all the
processes are continuous, stable convergence in law implies joint
convergence. Therefore considering \eqref{EqConvTriple}, we have
\[
\bigl(X,\ApproxVar{X}^{n},W,B^{2;n} \bigr)\Rightarrow
\bigl(X,X,W,B^{2} \bigr).
\]
Note that we proved tightness of the quadratic variation of the
martingale $\sqrt{n}B^{2;n}$, so that it has property (\ref{PropertyStar}), and
therefore it is good.

Now, $B^{3;n}$ is also a continuous Gaussian martingale, and we can
make use again of Theorem~4-1 in \citet{jacodcontinuous1997}. Let us
check the convergence in probability of its quadratic variation toward
$tI$ and of its quadratic covariation with respect to the other martingales.
Indeed, it is straightforward that if $j \neq j'$, $\langle
B^{3;i,j;n}, B^{3;i',j';n} \rangle_t =0$, while we deduce from the
multidimensional C.L.T. $\langle B^{3;i,j;n}, B^{3;i',j;n}\rangle_t
\mathop{\stackrel{P}{\longrightarrow}}\mathbf{1}_{\{i=i'\}}t.$ As
before this also shows goodness of $
B^{3;n}$. Using the same techniques we prove for any $i,j,l$ that
$\langle B^{3;l,j;n}, \sqrt{n} B^{2;i;n} \rangle_t =0, \mbox{ and
}\langle B^{3;l,j;n}, W \rangle=0$.
Hence
\[
\bigl(X,\ApproxVar{X}^{(n)},W,B^{2;n},B^{3;n}
\bigr)\Rightarrow \bigl(X,X,W,B^{2},B^{3} \bigr).
\]
We prove now the convergence in probability toward zero of the
remaining terms in the left side tuple in \eqref{EqTupleConv3}.

Since $n^{-1/2} \rightarrow0$ and $nA^{0,n} \Rightarrow A^{0}$, we
have $n^{1/2}A^{0,n}\Rightarrow0$ and thus $n^{1/2}\times  A^{0,n}\mathop
{\stackrel{P}{\longrightarrow}}0$.

On the other hand, for any $t\geq0$ and $k$,
\begin{eqnarray*}
{{\mathbb{E}}} \bigl[{ \bigl\langle\sqrt{n}B^{0;k;n}, \sqrt
{n}B^{0;k;n} \bigr\rangle_t }\bigr] & =& n \int
_0^t \bigl(s- \tbar[s] \bigr)^2 \,ds
\\
&=& \sum_{i=0}^{\lfloor nt\rfloor} n\int
_0^{1/n} r^2 \,dr + n\int
_0^{t-\lfloor nt\rfloor/n} r^2 \,dr
\\
& \leq&\sum_{i=0}^{\lfloor nt\rfloor+1} \frac{1}{3n^2}
= O \bigl(n^{-1} \bigr).
\end{eqnarray*}
So that by the Burckholder--Davis--Gundy inequality, ${{\mathbb{E}}}
[{\sup_{0\leq t \leq T} | \sqrt{n}B^{0;n} |^2 }]$ tends to zero as
$n\rightarrow
\infty$, implying $\sqrt{n}B^{0;n} \mathop{\stackrel
{P}{\longrightarrow}}0$.
In addition, it can be readily seen that
\[
{{\mathbb{E}}} \bigl[{\bigl|n A^{1;k;n} A^{1;j;n} \bigr|}\bigr] = 0
\]
for $j \neq k$, so that we have by using convex and Cauchy--Schwarz
inequalities,
\[
{{\mathbb{E}}} \Bigl[{\sup_{0\leq t \leq T}\bigl|\sqrt{n}A^{1;n}_t
\bigr|^2}\Bigr] \leq nT\sum_{j=1}^{d_x}
\int_0^T {{\mathbb{E}}} \bigl[{
\bigl(W_s^j-W_{\tbar{s}}^j \bigr)}
\bigr] \leq \frac{d_xT}{2},
\]
and hence, by the law of large numbers, $\sqrt{n}A^{1;n}\mathop
{\stackrel{P}{\longrightarrow}}{{\mathbb{E}}} [{\sqrt{n}A^{1;n}}] = 0$.

Finally, as $\sqrt{n}B^{1;n}$ converges in law to a Brownian,
$B^{3;n}\mathop{\stackrel{P}{\longrightarrow}}0$. Therefore \eqref
{EqTupleConv3} is proved.
\end{pf*}

\section{Cholesky decomposition}
\label{SecAnnexCholesky}

\begin{pf*}{Proof of Lemma~\ref{LemErrorCholesky}}
Since $G+\Delta G$ is the lower triangular factor of $H+\Delta H$, we have
\[
(G_{i,i} + \Delta G_{i,i})^2 = H_{i,i}
+ \Delta H_{i,i} - \sum_{k=1}^{i-1}(G_{i,k}
+\Delta G_{i,k})^2.
\]
By algebraic manipulation and the fact that $G$ is the Cholesky
decomposition of~$H$, we get
\[
\Delta G_{i,i} = \frac{\Delta H_{i,i} - 2 \sum_{k=1}^{i-1}\Delta
G_{i,k} G_{i,k}}{2G_{i,i}} - \Biggl( (\Delta
G_{i,i})^2 + \sum_{k=1}^{i-1}(
\Delta G_{i,k})^2 \Biggr).
\]
The first claim follows by controlling the last term by induction in
$i$, Theorem~\ref{TheoControlCholesky} and norm equivalence given by
\eqref{EqNormEquivalence}. The case $i>j$ is proved in the same way.
\end{pf*}
\end{appendix}

\section*{Acknowledgments}
The author would like to thank Fran\c{c}ois Delarue for his help and
support during the preparation of this work and the anonymous referee
for his suggestions that greatly improved the paper.
%
%



\printaddresses

\begin{thebibliography}{23}

\bibitem[\protect\citeauthoryear{Aldous and
Eagleson}{1978}]{aldousmixing1978}
%
\begin{barticle}[mr]
\bauthor{\bsnm{Aldous},~\bfnm{D.~J.}\binits{D.~J.}} \AND
\bauthor{\bsnm{Eagleson},~\bfnm{G.~K.}\binits{G.~K.}}
(\byear{1978}).
\btitle{On mixing and stability of limit theorems}.
\bjournal{Ann. Probab.}
\bvolume{6}
\bpages{325--331}.
\bid{mr={0517416}}
\end{barticle}
%
\bptok{imsref}%
\endbibitem

\bibitem[\protect\citeauthoryear{Bensoussan, Lions and
Papanicolaou}{1978}]{bensoussanasymptotic1978}
%
\begin{bbook}[mr]
\bauthor{\bsnm{Bensoussan},~\bfnm{Alain}\binits{A.}},
\bauthor{\bsnm{Lions},~\bfnm{Jacques-Louis}\binits{J.-L.}} \AND
\bauthor{\bsnm{Papanicolaou},~\bfnm{George}\binits{G.}}
(\byear{1978}).
\btitle{Asymptotic Analysis for Periodic Structures}.
\bseries{Studies in Mathematics and Its Applications}
\bvolume{5}.
\bpublisher{North-Holland},
\blocation{Amsterdam}.
\bid{mr={0503330}}
\end{bbook}
%
\bptok{imsref}%
\endbibitem

\bibitem[\protect\citeauthoryear{E, Liu and
Vanden-Eijnden}{2005}]{eanalysis2005}
%
\begin{barticle}[mr]
\bauthor{\bsnm{E},~\bfnm{Weinan}\binits{W.}},
\bauthor{\bsnm{Liu},~\bfnm{Di}\binits{D.}} \AND
\bauthor{\bsnm{Vanden-Eijnden},~\bfnm{Eric}\binits{E.}}
(\byear{2005}).
\btitle{Analysis of multiscale methods for stochastic differential equations}.
\bjournal{Comm. Pure Appl. Math.}
\bvolume{58}
\bpages{1544--1585}.
\bid{doi={10.1002/cpa.20088}, issn={0010-3640}, mr={2165382}}
\end{barticle}
%
\bptok{imsref}%
\endbibitem

\bibitem[\protect\citeauthoryear{E et~al.}{2007}]{eheterogeneous2007}
%
\begin{barticle}[mr]
\bauthor{\bsnm{E},~\bfnm{Weinan}\binits{W.}},
\bauthor{\bsnm{Engquist},~\bfnm{Bjorn}\binits{B.}},
\bauthor{\bsnm{Li},~\bfnm{Xiantao}\binits{X.}},
\bauthor{\bsnm{Ren},~\bfnm{Weiqing}\binits{W.}} \AND
\bauthor{\bsnm{Vanden-Eijnden},~\bfnm{Eric}\binits{E.}}
(\byear{2007}).
\btitle{Heterogeneous multiscale methods: A review}.
\bjournal{Commun. Comput. Phys.}
\bvolume{2}
\bpages{367--450}.
\bid{issn={1815-2406}, mr={2314852}}
\end{barticle}
%
\bptok{imsref}%
\endbibitem

\bibitem[\protect\citeauthoryear{Fouque and Lorig}{2011}]{fouquefast2011}
%
\begin{barticle}[mr]
\bauthor{\bsnm{Fouque},~\bfnm{Jean-Pierre}\binits{J.-P.}} \AND
\bauthor{\bsnm{Lorig},~\bfnm{Matthew~J.}\binits{M.~J.}}
(\byear{2011}).
\btitle{A fast mean-reverting correction to {H}eston's stochastic
volatility model}.
\bjournal{SIAM J. Financial Math.}
\bvolume{2}
\bpages{221--254}.
\bid{doi={10.1137/090761458}, issn={1945-497X}, mr={2775413}}
\end{barticle}
%
\bptok{imsref}%
\endbibitem

\bibitem[\protect\citeauthoryear{Fouque, Papanicolaou and
Sircar}{2000}]{fouquederivatives2000}
%
\begin{bbook}[mr]
\bauthor{\bsnm{Fouque},~\bfnm{Jean-Pierre}\binits{J.-P.}},
\bauthor{\bsnm{Papanicolaou},~\bfnm{George}\binits{G.}} \AND
\bauthor{\bsnm{Sircar},~\bfnm{K.~Ronnie}\binits{K.~R.}}
(\byear{2000}).
\btitle{Derivatives in Financial Markets with Stochastic Volatility}.
\bpublisher{Cambridge Univ. Press},
\blocation{Cambridge}.
\bid{mr={1768877}}
\end{bbook}
%
\bptok{imsref}%
\endbibitem

\bibitem[\protect\citeauthoryear{Fouque et~al.}{2003}]{fouquesingular2003}
%
\begin{barticle}[mr]
\bauthor{\bsnm{Fouque},~\bfnm{J.-P.}\binits{J.-P.}},
\bauthor{\bsnm{Papanicolaou},~\bfnm{G.}\binits{G.}},
\bauthor{\bsnm{Sircar},~\bfnm{R.}\binits{R.}} \AND
\bauthor{\bsnm{Solna},~\bfnm{K.}\binits{K.}}
(\byear{2003}).
\btitle{Singular perturbations in option pricing}.
\bjournal{SIAM J. Appl. Math.}
\bvolume{63}
\bpages{1648--1665}.
\bid{doi={10.1137/S0036139902401550}, issn={0036-1399}, mr={2001213}}
\end{barticle}
%
\bptok{imsref}%
\endbibitem

\bibitem[\protect\citeauthoryear{Jacod}{1997}]{jacodcontinuous1997}
%
\begin{bincollection}[mr]
\bauthor{\bsnm{Jacod},~\bfnm{Jean}\binits{J.}}
(\byear{1997}).
\btitle{On continuous conditional {G}aussian martingales and stable
convergence in law}.
In \bbooktitle{S\'eminaire de {P}robabilit\'es, {XXXI}}.
\bseries{Lecture Notes in Math.}
\bvolume{1655}
\bpages{232--246}.
\bpublisher{Springer},
\blocation{Berlin}.
\bid{doi={10.1007/BFb0119308}, mr={1478732}}
\end{bincollection}
%
\bptok{imsref}%
\endbibitem

\bibitem[\protect\citeauthoryear{Jacod and
Protter}{1998}]{jacodasymptotic1998}
%
\begin{barticle}[mr]
\bauthor{\bsnm{Jacod},~\bfnm{Jean}\binits{J.}} \AND
\bauthor{\bsnm{Protter},~\bfnm{Philip}\binits{P.}}
(\byear{1998}).
\btitle{Asymptotic error distributions for the {E}uler method for
stochastic differential equations}.
\bjournal{Ann. Probab.}
\bvolume{26}
\bpages{267--307}.
\bid{doi={10.1214/aop/1022855419}, issn={0091-1798}, mr={1617049}}
\end{barticle}
%
\bptok{imsref}%
\endbibitem

\bibitem[\protect\citeauthoryear{Jakubowski, M{\'e}min and Pag{\`
e}s}{1989}]{jakubowskiconvergence1989}
%
\begin{barticle}[mr]
\bauthor{\bsnm{Jakubowski},~\bfnm{A.}\binits{A.}},
\bauthor{\bsnm{M{\'e}min},~\bfnm{J.}\binits{J.}} \AND
\bauthor{\bsnm{Pag{\`e}s},~\bfnm{G.}\binits{G.}}
(\byear{1989}).
\btitle{Convergence en loi des suites d'int\'egrales stochastiques
sur l'espace {$\mathbf{D}^1$} de {S}korokhod}.
\bjournal{Probab. Theory Related Fields}
\bvolume{81}
\bpages{111--137}.
\bid{doi={10.1007/BF00343739}, issn={0178-8051}, mr={0981569}}
\end{barticle}
%
\bptok{imsref}%
\endbibitem

\bibitem[\protect\citeauthoryear{Kurtz and Protter}{1991a}]{kurtzweak1991}
%
\begin{barticle}[mr]
\bauthor{\bsnm{Kurtz},~\bfnm{Thomas~G.}\binits{T.~G.}} \AND
\bauthor{\bsnm{Protter},~\bfnm{Philip}\binits{P.}}
(\byear{1991}a).
\btitle{Weak limit theorems for stochastic integrals and stochastic
differential equations}.
\bjournal{Ann. Probab.}
\bvolume{19}
\bpages{1035--1070}.
\bid{issn={0091-1798}, mr={1112406}}
\end{barticle}
%
\bptok{imsref}%
\endbibitem

\bibitem[\protect\citeauthoryear{Kurtz and
Protter}{1991b}]{kurtzwong-zakai1991}
%
\begin{bincollection}[mr]
\bauthor{\bsnm{Kurtz},~\bfnm{Thomas~G.}\binits{T.~G.}} \AND
\bauthor{\bsnm{Protter},~\bfnm{Philip}\binits{P.}}
(\byear{1991}b).
\btitle{Wong--{Z}akai corrections, random evolutions, and simulation
schemes for {SDE}s}.
In \bbooktitle{Stochastic Analysis}
\bpages{331--346}.
\bpublisher{Academic Press},
\blocation{Boston, MA}.
\bid{mr={1119837}}
\end{bincollection}
%
\bptok{imsref}%
\endbibitem

\bibitem[\protect\citeauthoryear{Kurtz and Protter}{1996}]{kurtzweak1996}
%
\begin{bincollection}[mr]
\bauthor{\bsnm{Kurtz},~\bfnm{Thomas~G.}\binits{T.~G.}} \AND
\bauthor{\bsnm{Protter},~\bfnm{Philip~E.}\binits{P.~E.}}
(\byear{1996}).
\btitle{Weak convergence of stochastic integrals and differential equations}.
In \bbooktitle{Probabilistic Models for Nonlinear Partial Differential
Equations ({M}ontecatini {T}erme, 1995)}.
\bseries{Lecture Notes in Math.}
\bvolume{1627}
\bpages{1--41}.
\bpublisher{Springer},
\blocation{Berlin}.
\bid{doi={10.1007/BFb0093176}, mr={1431298}}
\end{bincollection}
%
\bptok{imsref}%
\endbibitem

\bibitem[\protect\citeauthoryear{Lamberton and Pag{\`
e}s}{2002}]{lambertonrecursive2002}
%
\begin{barticle}[mr]
\bauthor{\bsnm{Lamberton},~\bfnm{Damien}\binits{D.}} \AND
\bauthor{\bsnm{Pag{\`e}s},~\bfnm{Gilles}\binits{G.}}
(\byear{2002}).
\btitle{Recursive computation of the invariant distribution of a diffusion}.
\bjournal{Bernoulli}
\bvolume{8}
\bpages{367--405}.
\bid{doi={10.1142/S0219493703000838}, issn={1350-7265}, mr={1913112}}
\end{barticle}
%
\bptok{imsref}%
\endbibitem

\bibitem[\protect\citeauthoryear{Lemaire}{2005}]{lemaireestimation2005}
%
\begin{bmisc}[author]
\bauthor{\bsnm{Lemaire},~\bfnm{Vincent}\binits{V.}}
(\byear{2005}).
\bhowpublished{Estimation r{\'e}cursive de la mesure invariante d'un
processus de diffusion.
Ph.D. thesis,
Univ. de Marne-la-Vall{\'e}e}.
\end{bmisc}
%
\bptok{imsref}%
\endbibitem

\bibitem[\protect\citeauthoryear{Majda, Timofeyev and
Vanden-Eijnden}{2001}]{majdamathematical2001}
%
\begin{barticle}[mr]
\bauthor{\bsnm{Majda},~\bfnm{Andrew~J.}\binits{A.~J.}},
\bauthor{\bsnm{Timofeyev},~\bfnm{Ilya}\binits{I.}} \AND
\bauthor{\bsnm{Vanden-Eijnden},~\bfnm{Eric}\binits{E.}}
(\byear{2001}).
\btitle{A mathematical framework for stochastic climate models}.
\bjournal{Comm. Pure Appl. Math.}
\bvolume{54}
\bpages{891--974}.
\bid{doi={10.1002/cpa.1014}, issn={0010-3640}, mr={1829529}}
\end{barticle}
%
\bptok{imsref}%
\endbibitem

\bibitem[\protect\citeauthoryear{Meyn and Tweedie}{1993}]{meynstability1993}
%
\begin{barticle}[mr]
\bauthor{\bsnm{Meyn},~\bfnm{Sean~P.}\binits{S.~P.}} \AND
\bauthor{\bsnm{Tweedie},~\bfnm{R.~L.}\binits{R.~L.}}
(\byear{1993}).
\btitle{Stability of {M}arkovian processes. {III}.
{F}oster--{L}yapunov criteria for continuous-time processes}.
\bjournal{Adv. in Appl. Probab.}
\bvolume{25}
\bpages{518--548}.
\bid{doi={10.2307/1427522}, issn={0001-8678}, mr={1234295}}
\end{barticle}
%
\bptok{imsref}%
\endbibitem

\bibitem[\protect\citeauthoryear{Pardoux and
Veretennikov}{2001}]{pardouxpoisson2001}
%
\begin{barticle}[mr]
\bauthor{\bsnm{Pardoux},~\bfnm{E.}\binits{E.}} \AND
\bauthor{\bsnm{Veretennikov},~\bfnm{A.~Y.}\binits{A.~Y.}}
(\byear{2001}).
\btitle{On the {P}oisson equation and diffusion approximation. {I}}.
\bjournal{Ann. Probab.}
\bvolume{29}
\bpages{1061--1085}.
\bid{doi={10.1214/aop/1015345596}, issn={0091-1798}, mr={1872736}}
\end{barticle}
%
\bptok{imsref}%
\endbibitem

\bibitem[\protect\citeauthoryear{Pardoux and
Veretennikov}{2003}]{pardouxpoisson2003}
%
\begin{barticle}[mr]
\bauthor{\bsnm{Pardoux},~\bfnm{{\`E}.}\binits{{\`E}.}} \AND
\bauthor{\bsnm{Veretennikov},~\bfnm{A.~Y.}\binits{A.~Y.}}
(\byear{2003}).
\btitle{On {P}oisson equation and diffusion approximation. {II}}.
\bjournal{Ann. Probab.}
\bvolume{31}
\bpages{1166--1192}.
\bid{doi={10.1214/aop/1055425774}, issn={0091-1798}, mr={1988467}}
\end{barticle}
%
\bptok{imsref}%
\endbibitem

\bibitem[\protect\citeauthoryear{Sun}{1991}]{sunperturbation1991}
%
\begin{barticle}[mr]
\bauthor{\bsnm{Sun},~\bfnm{Ji~Guang}\binits{J.~G.}}
(\byear{1991}).
\btitle{Perturbation bounds for the {C}holesky and {$QR$} factorizations}.
\bjournal{BIT}
\bvolume{31}
\bpages{341--352}.
\bid{doi={10.1007/BF01931293}, issn={0006-3835}, mr={1112229}}
\end{barticle}
%
\bptok{imsref}%
\endbibitem

\bibitem[\protect\citeauthoryear{Talay}{1990}]{Talaysecond-order1990}
%
\begin{barticle}[author]
\bauthor{\bsnm{Talay},~\bfnm{Denis}\binits{D.}}
(\byear{1990}).
\btitle{Second-order discretization schemes of stochastic
differential systems for the computation of the invariant law}.
\bjournal{Stochastics and Stochastic Reports}
\bvolume{29}
\bpages{13--36}.
\end{barticle}
%
\bptok{imsref}%
\endbibitem

\bibitem[\protect\citeauthoryear
{Veretennikov}{1997}]{veretennikovpolynomial1997}
%
\begin{barticle}[mr]
\bauthor{\bsnm{Veretennikov},~\bfnm{A.~Y.}\binits{A.~Y.}}
(\byear{1997}).
\btitle{On polynomial mixing bounds for stochastic differential equations}.
\bjournal{Stochastic Process. Appl.}
\bvolume{70}
\bpages{115--127}.
\bid{doi={10.1016/S0304-4149(97)00056-2}, issn={0304-4149}, mr={1472961}}
\end{barticle}
%
\bptok{imsref}%
\endbibitem

\bibitem[\protect\citeauthoryear
{Veretennikov}{2011}]{veretennikovsobolev2011}
%
\begin{barticle}[mr]
\bauthor{\bsnm{Veretennikov},~\bfnm{A.~Y.}\binits{A.~Y.}}
(\byear{2011}).
\btitle{On {S}obolev solutions of {P}oisson equations in {$\mathbb
{R}^d$} with a parameter}.
\bjournal{J. Math. Sci. (N.Y.)}
\bvolume{179}
\bpages{48--79}.
\bid{doi={10.1007/s10958-011-0582-5}, issn={1072-3374}, mr={3014098}}
\end{barticle}
%
\bptok{imsref}%
\endbibitem
\end{thebibliography}
\end{document}